\documentclass{amsart}
\input epsf

\usepackage{amsfonts,amsthm,amsmath,amssymb,latexsym}
\usepackage{graphicx,color}
\usepackage[all]{xy}

\begin{document}

\author{Dragomir \v Sari\' c}
\thanks{This research is partially supported by National Science Foundation grant DMS 1102440.}

\address{Department of Mathematics, Queens College of CUNY,
65-30 Kissena Blvd., Flushing, NY 11367}
\email{Dragomir.Saric@qc.cuny.edu}

\address{Mathematics PhD. Program, The CUNY Graduate Center, 365 Fifth Avenue, New York, NY 10016-4309}

\theoremstyle{definition}

 \newtheorem{definition}{Definition}[section]
 \newtheorem{remark}[definition]{Remark}
 \newtheorem{example}[definition]{Example}

\newtheorem*{notation}{Notation}

\theoremstyle{plain}

 \newtheorem{proposition}[definition]{Proposition}
 \newtheorem{theorem}[definition]{Theorem}
 \newtheorem{corollary}[definition]{Corollary}
 \newtheorem{lemma}[definition]{Lemma}

\def\H{{\mathbf H}}
\def\F{{\mathcal F}}
\def\R{{\mathbf R}}
\def\Q{{\mathbf Q}}
\def\Z{{\mathbf Z}}
\def\E{{\mathcal E}}
\def\N{{\mathbf N}}
\def\X{{\mathcal X}}
\def\Y{{\mathcal Y}}
\def\C{{\mathbf C}}
\def\D{{\mathbf D}}
\def\G{{\mathcal G}}

\title{Complex Fenchel-Nielsen coordinates with small imaginary parts}


\date{\today}

\maketitle

\begin{abstract}
Kahn and Markovic \cite{KahnMark} proved that the fundamental group of each closed hyperbolic three manifold contains a closed surface subgroup. One of the main ingredients in their proof is a theorem which states that an assignment of nearly real, complex Fenchel-Nielsen coordinates to the cuffs of a pants decomposition of a closed surface $S$ induces a quasiFuchsian representation of the fundamental group of $S$. We give a new proof of this theorem with a slightly stronger conditions on the Fenchel-Nielsen coordinates and explain how to use the exponential mixing of the geodesic flow on a closed hyperbolic three manifold to prove that our theorem is sufficient for the applications in the work of Kahn and Markovic \cite{KahnMark}. 
\end{abstract}

\section{Introduction}

Kahn and Markovic \cite{KahnMark} recently proved that the fundamental group $\pi_1(M)$ of any closed hyperbolic three
manifold $M$ has a closed surface subgroup. Their  proof uses the exponential 
mixing of the geodesic flow on $M$ in order to find
a ``well-distributed'' finite collection of skew pants in the three manifold $M$ that have large and nearly equal cuff lengths, that are nearly flat, and that can be glued pairwise with nearly zero angles.
The collection
of skew pants has a subcollection $\mathcal{P}$ that closes to form an abstract closed surface $S$ of (high) genus with nearly real, complex Fenchel-Nielsen coordinates on $\mathcal{P}$. 
The final step in the proof of Kahn and Markovic \cite{KahnMark} is
to show that nearly real, complex Fenchel-Nielsen coordinates on $\mathcal{P}$ necessarily induce an isomorphism between the fundamental group $\pi_1(S)$ and a quasiFuchsian group. Our contribution is to give a new proof of this statement. In fact, we prove a slightly weaker statement by requiring that the imaginary parts of the complex Fenchel-Nielsen coordinates satisfy slightly stronger conditions and establish that this weaker statement is sufficient  for the purposes of the proof of the surface subgroup conjecture along the lines in \cite{KahnMark}. Our proof adopts the ideas of proving the injectivity of the bending along a measured lamination (cf. \cite{EpMarMar}, \cite{KeeSe}, \cite{Sa1}) which (at least conceptually) simplifies this part of the argument in \cite{KahnMark}.

Let $S$ be a closed surface of genus $g\geq 2$ equipped with a pants decomposition $\mathcal{P}$.  Then $\mathcal{P}$ consists of $3g-3$ simple closed curves such that each component of the complement is a pair of pants. 
Following \cite[\S 2]{KahnMark} (see also \S 2), to each cuff $C\in\mathcal{P}$ we associate complex half-length $hl(C)\in\mathbb{C}/2\pi i \mathbb{Z}$ and complex twist-bend parameter $s(C)\in\mathbb{C}/(2\pi i\mathbb{Z}+hl(C)\mathbb{Z})$. 
An assignment of half-lengths and twist-bend parameters to $C\in\mathcal{P}$ induces a representation of the fundamental group $\pi_1(S)$ into $PSL_2(\mathbb{C})$. The representation is Fuchsian if and only if $\{ (hl(C),s(C)\}_{C\in\mathcal{P}}\in\mathbb{R}^{3g-3}$. The following theorem characterizes a neighborhood of the real subspace $\mathbb{R}^{3g-3}$ inside $\mathbb{C}^{3g-3}$ which gives quasiFuchsian representations.
We note that the size of the neighborhood is independent of the genus $g$.

\begin{theorem}[Kahn-Markovic \cite{KahnMark}]
\label{thm:mark_kahn}
There exist universal $\hat{\epsilon}, K_0>0$ and $R(\hat{\epsilon})>0$ such that the following is satisfied. Let $S$ be a closed surface of genus $g\geq 2$ and let $\mathcal{P}$ be a pants decomposition of $S$.
If $\rho :\pi_1(S)\to PSL_2(\mathbb{C})$ is a representation which is discrete and faithful on each pair of pants in $\mathcal{P}$ and if the reduced complex Fenchel-Nielsen coordinates on each cuff $C\in\mathcal{P}$
satisfy
\begin{equation}
\label{eq:half-length}
|hl(C)-R/2|<\epsilon
\end{equation}
and
\begin{equation}
\label{eq:twist-bend}
|s(C)-1|<\epsilon /R
\end{equation}
for some $\epsilon<\hat{\epsilon}$ and $R>R(\hat{\epsilon})$, then $\rho :\pi_1(S\to PSL_2(\mathbb{C})$ is injective and $\rho (\pi_1(S))$ is quasiFuchsian.

Moreover, let $S$ be endowed with a hyperbolic metric whose reduced Fenchel-Nielsen coordinates are $hl(C)=R/2$ and $s(C)=1$ for each $C\in\mathcal{P}$. Then there exists an injective map $\tilde{f}:\partial_{\infty}\tilde{S}\to\partial_{\infty}\mathbb{H}^3$ which conjugates $\pi_1(S)$ into the above quasiFuchsian group and which extends to a $(1+K_0\epsilon )$-quasiconformal map of $\partial_{\infty}\mathbb{H}^3$ onto itself.
\end{theorem}

We give a new proof of the above theorem under assumptions (\ref{eq:twist-bend}) and
\begin{equation}
\label{eq:half-length-modif}
|hl(C)-R/2|<\epsilon /R.
\end{equation}

Even though (\ref{eq:half-length-modif}) is stronger than (\ref{eq:half-length}) (which makes our statement weaker than the above theorem), it turns out that this is enough for the purposes in \cite{KahnMark}. At the end of Introduction we indicate how to see that 
(\ref{eq:half-length-modif}) follows from the fact that the skew pairs of pants are ``well-distributed'' inside the three manifold which proves that the weaker statement suffices. One advantage of using  (\ref{eq:half-length-modif}) instead of (\ref{eq:half-length}) is that we do not need to require that $\rho$ is discrete and faithful on pairs of pants in $\mathcal{P}$ in order to establish the injectivity of the representation $\rho :\pi_1(S)\to PSL_2(\mathbb{C})$. In particular, $\rho$ is discrete and injective on each pair of pants of $\mathcal{P}$ if it satisfies (\ref{eq:twist-bend}) and (\ref{eq:half-length-modif}).

In \cite{KahnMark}, Theorem \ref{thm:mark_kahn} is proved by estimating the derivative (along a path of representations connecting the Fuchsian representation with the representation corresponding to the given reduced Fenchel-Nielsen coordinates) of the distance between the images in $\mathbb{H}^3$ of two lifts of geodesics $C\in\mathcal{P}$ in $\mathbb{H}^2$ from the above by a function of the distance between these two lifts of geodesics at the representation.  
This leads to an inductive argument which gives the desired theorem. 

Our approach is to decompose each pair of pants in $\mathcal{P}$ into two ideal hyperbolic triangles by adding three infinite geodesics such that each end of each added geodesic spirals around a different cuff. The union of cuffs of $\mathcal{P}$ together with the added geodesics in each pair of pants is a maximal geodesic lamination $\lambda$ in $S$ with finitely many leaves. Let $\tilde{\lambda}$ be the lift of $\lambda$ to the universal covering $\mathbb{H}^2$. The reduced complex Fenchel-Nielsen coordinates $\{ (hl(C),s(C)\}_{C\in\mathcal{P}}$ induce a developing map $\tilde{f}:\partial_{\infty}\mathbb{H}^2\to\partial_{\infty}\mathbb{H}^3$ which conjugates $\pi_1(S)<PSL_2(\mathbb{R})$ into a subgroup of $PSL_2(\mathbb{C})$. The developing map extends to complementary triangles of $\tilde{\lambda}$ to define a pleated surface $\tilde{f}: \mathbb{H}^2\to \mathbb{H}^3$ which is pleated along $\tilde{\lambda}$ (cf. \cite{Bon2} and also \S 3). Each pleated surface along $\tilde{\lambda}$ induces a finitely additive $(\mathbb{C}/2\pi i\mathbb{Z})$-valued transverse cocycle $\alpha$ to $\tilde{\lambda}$ which measures the shearing and the bending along $\tilde{\lambda}$ (cf. \cite{Bon2}). The {\it bending cocycle} is the imaginary part $\beta$ of the transverse cocycle $\alpha$ which is  an $(\mathbb{R}/2\pi\mathbb{Z})$-valued transverse cocycle to $\tilde{\lambda}$ measuring the amount of the bending of the pleated surface. We translate the reduced complex Fenchel-Nielsen coordinates into the bending transverse cocycle $\beta$ to $\tilde{\lambda}$ as follows. 

An isolated leaf $\tilde{l}$ of $\tilde{\lambda}$ is on a common boundary of two complementary ideal triangles $\Delta^{1}(\tilde{l})$ and $\Delta^2(\tilde{l})$ to $\tilde{\lambda}$. Isolated leaves of $\tilde{\lambda}$ accumulate to each lift $\tilde{C}$ of each cuff $C\in\mathcal{P}$ from both sides of $\tilde{C}$. Let $S$ be endowed with a hyperbolic metric whose Fenchel-Nielsen coordinates are $\{ (Re(hl(C)),Re(s(C))\}_{C\in\mathcal{P}}$ and divide each pair of pants of $\mathcal{P}$ into two hyperbolic hexagons $\Sigma_1$ and $\Sigma_2$ by drawing common orthogonal arcs between pairs of cuffs of each pair of pants in $\mathcal{P}$.  Each hexagon $\Sigma$ on the surface $S$ lifts to infinitely many hexagons in the universal covering $\mathbb{H}^2$. 
Recall the assumptions $|Re(s(C))-1|<\epsilon /R$ and $|Re(hl(C))-R/2|<\epsilon /R$, for $R\geq R(\hat{\epsilon})$ and $0<\epsilon <\hat{\epsilon}$. A cuff $C$ is the union of two boundary sides of two hexagons coming from the pair of pants on one side of $C$ as well as the union of two boundary sides of two hexagons coming from the pair of pants on the other side of $C$. The boundary sides of the hexagons from one side of $C$ are not exactly matched along $C$ with the boundaries of the hexagons from the other side of $C$ but they are glued with a shift close to $1$ by the condition $|Re(s(C))-1|<\epsilon /R$. It follows that for each hexagon $\Sigma_i$, $i=1,2$, on one side of $C$ there is a unique hexagon $\Sigma_i'$, $i=1,2$, on the other side of $C$ such that the common subarc of their boundary sides on $C$ has length close to $R/2-1$.  We will say that $\Sigma_i$ and $\Sigma_i'$ are $0$-neighbors in this case. Two lifts $\tilde{\Sigma}_i$ and $\tilde{\Sigma}_i'$ to $\mathbb{H}^2$ of $0$-neighbors hexagons are also called $0$-neighbors if they meet along a lift $\tilde{C}$ of $C$ with a common subarc of length close to $R/2-1$ (cf. \S 3).  If hexagon $\tilde{\Sigma}$ is a lift of a hexagon $\Sigma$, then 
$\tilde{\Sigma}$ intersects infinitely many complementary triangles to $\tilde{\lambda}$. There is a unique triangle $\Delta_{\tilde{\Sigma}}$ such that its intersection with $\tilde{\Sigma}$ is a hexagon, and we call $\Delta_{\tilde{\Sigma}}$ the {\it canonical triangle} of $\tilde{\Sigma}$ (cf. \S 3 and Figure 1).

\begin{theorem}
\label{thm:f-n into shear-bend} There exists $C_0>0$ such that the following holds. Let $\{
hl(C),s(C)\}_{C\in\mathcal{P}}$ be the reduced complex Fenchel-Nielsen
coordinates that satisfy {\rm (\ref{eq:half-length-modif})} and {\rm (\ref{eq:twist-bend})}, and let $\beta$ be the induced bending transverse cocycle to the lamination $\tilde{\lambda}$.  If $\tilde{l}$ is an isolated leaf of $\tilde{\lambda}$ and $\Delta^{i}(\tilde{l})$, $i=1,2$, are complementary triangles to $\tilde{\lambda}$ with a common boundary side $\tilde{l}$, then
\begin{equation}
\label{eq:bending-isolated-leaf}
|\beta (\Delta^1(\tilde{l}),\Delta^2(\tilde{l}))|\leq \frac{C_0\epsilon}{R} .
\end{equation}
If $\tilde{\Sigma}_1$ and $\tilde{\Sigma}_2$ are $0$-neighbors hexagons, and $\Delta_{\tilde{\Sigma}_1}$ and $\Delta_{\tilde{\Sigma}_2}$ are their canonical
triangles, then
\begin{equation}
\label{eq:bending-cuff}
|\beta (\Delta_{\tilde{\Sigma}_1},\Delta_{\tilde{\Sigma}_2})|\leq
\frac{C_0\epsilon}{R}.
\end{equation}
\end{theorem}

Theorem \ref{thm:f-n into shear-bend} translates the original problem of whether the reduced complex Fenchel-Nielsen coordinates give a quasiFuchsian representation into whether a bending cocycle gives a quasiFuchsian representation. We point out that any condition on the bending cocycle that guarantees injectivity of the bending map on $\partial_{\infty}\mathbb{H}^2$ necessarily depends on the hyperbolic metric from which the bending starts. In our case, a sufficient information about the hyperbolic metric is given by the fact that a unit length hyperbolic arc in $\mathbb{H}^2$ can have at most $2R+2$ intersections with the lifts of the cuffs (cf. \cite[Lemma 2.3]{KahnMark} and Lemma \ref{lem:length1-intersects-at-most-R}).  When the transverse bending measure is countably additive, then it is well-known that the bending measure which is uniformly small on each transverse arc of length $1$ gives a bending map which is injective on $\partial_{\infty}\mathbb{H}^2$ (cf. \cite{EpMarMar}, \cite{KeeSe}, \cite{Sa1}). The above sufficient condition for the injectivity of the bending along the measured laminations is universal in the sense that it does not depend on the genus of the surface and, in fact, it works for any surface including the unit disk.  
The main difficulty in proving that the bending along finitely additive transverse cocycles is injective on $\partial_{\infty}\mathbb{H}^2$ lies in the fact that the total variation of the transverse bending measure is infinite. However, conditions (\ref{eq:half-length-modif}) and (\ref{eq:twist-bend}) guarantee that the bending map behaves ``semi-locally'' as the bending map along a measured lamination. 
Using ideas from \cite{Sa1}, we prove that the above conditions on the bending cocycle and the hyperbolic metric are sufficient to guarantee injectivity on $\partial_{\infty}\mathbb{H}^2$ of the bending map. Holomorphic motions provide the desired bound on quasiconformal extension of the restriction to $\partial_{\infty}\mathbb{H}^2$ of the bending map.

\begin{theorem}
\label{thm:injective-bending-cocycles}
Given $C_0$, there exist $\hat{\epsilon}>0$, $K_0> 1$ and $R(\hat{\epsilon})>0$ such that for $0\leq\epsilon <\hat{\epsilon}$ the following is satisfied. Let $S$ be a closed hyperbolic surface equipped with a maximal, finite geodesic lamination $\lambda$ such that each closed geodesic of $\lambda$ has length in the interval $(R-\frac{C_0\epsilon}{R},R+\frac{C_0\epsilon}{R})$ for some $R\geq R(\hat{\epsilon})$ and that each geodesic arc on $S$ of length $1$ intersects at most $C_0\cdot R$ closed geodesics of $\lambda$.  If a bending cocycle $\beta$ transverse to the lift $\tilde{\lambda}$ in $\mathbb{H}^2$ satisfies (\ref{eq:bending-isolated-leaf}) and (\ref{eq:bending-cuff})
then the induced bending map
$$
\tilde{f}_{\beta}:\partial_{\infty}\mathbb{H}^2\to\partial_{\infty}\mathbb{H}^3
$$
is injective and the induced representation of $\pi_1(S)$ is quasiFuchsian. The bending map extends to a $(1+K_0\epsilon)$-quasiconformal map $\tilde{f}_{\beta}:\partial_{\infty}\mathbb{H}^3\to \partial_{\infty}\mathbb{H}^3$. 
\end{theorem} 

We give an analogue of the above theorem for non-finite geodesic laminations and bending (finitely additive) transverse cocycles \cite{Sa2}. 

\vskip .2 cm

It remains to explain why the condition $|hl(C)-R/2|<\epsilon$ can be replaced with the condition $|hl(C)-R/2|<\epsilon /R$. The geodesic flow on a closed hyperbolic three manifold $M$ is exponentially mixing \cite{Moo}. Let $\mathcal{F}(M)$ be the $2$-frame bundle over $M$. Since the hyperbolic Laplacian of $M$ has a spectral gap the following holds by \cite{Moo}.
There exists $q>0$ which depends on the three manifold $M$ such that for any two $C^{nifty}$-functions $\psi ,\varphi :\mathcal{F}(M)\to\mathbb{R}$ we have
\begin{equation}
\label{eq:expo}
\Big{|}\Lambda(\mathcal{F}(M))\int_{\mathcal{F}(M)}(g^{*}_t\psi)\varphi d\Lambda-\int_{\mathcal{F}(M)}\psi d\Lambda \int_{\mathcal{F}(M)}\varphi d\Lambda\Big{|}\leq Ce^{-qt},
\end{equation}
where $\Lambda$ is the Liouville measure on $\mathcal{F}$ and the constant $C$ depends on the $C^1$-norms of $\psi$ and $\varphi$. 

Let $f_{\epsilon}:\mathcal{F}(M)\to\mathbb{R}$ be a non-negative $C^{\infty}$-function 
supported in the
$\epsilon$-neighborhood of a point in $\mathcal{F}(M^3)$ with $\int_{\mathcal{F}(M)}f_{\epsilon}d\Lambda =1$, called a {\it bump function} for the $\epsilon$-neighborhood. By applying (\ref{eq:expo}) to $f_{\epsilon}$, Kahn and Markovic \cite{KahnMark} proved that there exist triples of $2$-frames, called {\it tripods}, that after traveling a long time $t>0$ along the geodesic flow return to their $\epsilon$-neighborhoods. These tripods define skew pairs of pants in $M$ whose cuff lengths are $R=2t-2\log \frac{4}{3}$ with a possible error $D\cdot \epsilon$, for some $D>0$  because the expression on the right of (\ref{eq:expo}) goes to $0$ as $t\to\infty$ (cf. \cite[Lemma 4.6]{KahnMark}). 

It is possible to improve the estimate on the complex length of the cuffs of the above skew pairs of pants. Note that the constant $C$ in (\ref{eq:expo}) can be estimated (cf. \cite{Moo}) in terms of $H^2_2$-Sobolev norms to be $C_0\|\psi\|_{H^2_2} \|\varphi\|_{H^2_2}$, where $C_0$ is a fixed constant.
Then, for a given time $t$ geodesic flow, we consider bump function $f_{\epsilon /t}$ for $\epsilon /t$-neighborhood of a point in $\mathcal{F}(M)$. The bump function $f_{\epsilon /t}$ can be produced by scalling the domain and the size of $f_{\epsilon}$ such that $\| f_{\epsilon /t}\|_{H^2_2}\leq p_2(t) \| f_{\epsilon}\|_{H^2_2}$, where $p_2(t)$ is a polynomial in $t$ of degree $2$. From (\ref{eq:expo}) we get
$$
\Big{|}\Lambda(\mathcal{F}(M))\int_{\mathcal{F}(M)}(g^{*}_tf_{\epsilon /t})f_{\epsilon /t} d\Lambda-1\Big{|}\leq C_0\| f_{\epsilon}\|_{H^2_2}^2p_2(t)^2e^{-qt}\to 0
$$
as $t\to\infty$. This implies that the skew pairs of pants have cuffs of the length $R$ within $\epsilon /R$, when we choose an appropriate value for $t=t(R)$ thus obtaining (\ref{eq:half-length-modif}). 

\vskip .2 cm

\paragraph{\it Acknowledgements} I am grateful to Vladimir Markovic and Hideki Miyachi for various discussions regarding this work. 

\section{The reduced complex Fenchel-Nielsen coordinates}

Let $S$ be a closed surface of genus $g\geq 2$ and let $\pi_1(S)$ be its fundamental group. Let $\mathcal{P}$ be a pants decomposition of $S$, namely $\mathcal{P}$ consists of $3g-3$ simple, closed curves on $S$ such that the components of the complement of $\mathcal{P}$ are pairs of pants. A representation $\rho :\pi_1(S)\to PSL_2(\mathbb{C})$ associates to each cuff $C\in\mathcal{P}$ two complex numbers: the complex length and the twist-bend parameter. In total, $6g-6$ complex numbers is associated to a representation $\rho$, called the complex Fenchel-Nielsen coordinates.
The complex Fenchel-Nielsen coordinates were introduced in \cite{Kou} and \cite{tan}, and it was proved there that the quasiFuchsian space of $S$ is parametrized by an open subset of $\mathbb{C}^{6g-6}$ which contains the real locus $\mathbb{R}^{6g-6}$. We use the reduced Fenchel-Nielsen coordinates introduced by Kahn and Markovic \cite[\S 2]{KahnMark} and we refer the reader to their article for more details.

Let $\alpha$ and $\beta$ be two oriented geodesics in $\mathbb{H}^3$. Let $\gamma$ be their common orthogonal oriented from $\alpha$ to $\beta$. The complex distance $d_{\gamma}(\alpha ,\beta )$ between $\alpha$ and $\beta$ is defined to have a positive real part equal to the distance between $\alpha\cap\gamma$ and $\beta\cap\gamma$, while the imaginary part of  $d_{\gamma}(\alpha ,\beta )$ is the angle between the parallel transport along $\gamma$ of the unit tangent vector to $\alpha$ at $\alpha\cap \gamma$ and the unit tangent vector to $\beta$ at $\beta\cap\gamma$. Since the imaginary part of $d_{\gamma}(\alpha ,\beta )$ is well defined modulo $2\pi i$, we have $d_{\gamma}(\alpha ,\beta )\in\mathbb{C}/2\pi i\mathbb{Z}$ (for more details, cf. \cite[\S 2]{KahnMark}).

Let $\Pi_0$ be a pair of pants and $\pi_1(\Pi_0)$ be its fundamental group. Consider a representation $\rho :\pi_1(\Pi_0)\to PSL_2(\mathbb{C})$ which is faithful and loxodromic. Namely, the cuffs $C_i$, $i=0,1,2$, of $\Pi_0$ are represented by loxodromic elements $\rho (C_i)\in PSL_2(\mathbb{C})$. Let $\gamma_i$ be the axis of $\rho (C_i)$ and $\eta_i$ be the common orthogonal to $\gamma_{i-1}$ and $\gamma_{i+1}$, for $i=0,1,2$, where the indices are taken modulo $3$. Then the half-length $hl_{\Pi_0 ,\rho}(C_i)$ of the curve $C_i$ associated to the representation $\rho$ is $d_{\gamma_i}(\eta_{i-1},\eta_{i+1})$ (cf. \cite[\S 2]{KahnMark}).

Consider a representation $\rho :\pi_1(S)\to PSL_2(\mathbb{C})$ of the
fundamental group $\pi_1(S)$ of a closed surface $S$ of genus at
least two into $PSL_2(\mathbb{C})$ and fix a pants of decomposition $\mathcal{P}$. The representation $\rho$ is {\it viable} if $\rho :\pi_1(\Pi )\to PSL_2(\mathbb{C})$ is discrete and faithful for each pair of pants $\Pi$ of the pants decomposition $\mathcal{P}$, and for any two pairs of pants $\Pi$ and $\Pi'$ with a common cuff $C$ we have $hl_{\Pi ,\rho}(C)=hl_{\Pi',\rho}(C)$ (cf. \cite{KahnMark}). For a given viable representation $\rho$, we define the complex half-length of a cuff $C\in\mathcal{P}$ by $hl(C)=hl_{\Pi ,\rho}(C)=hl_{\Pi',\rho}(C)\in\mathbb{C}/2\pi i\mathbb{Z}$, where $\Pi$ and $\Pi'$ are pairs of pants with one cuff $C$ (cf. \cite{KahnMark}). 
Let $\Pi$ and $\Pi'$ be two pairs of pants with cuffs $C_i$, $i=0,1,2$, and $C_i'$, $i=0,1,2$, respectively such that $C_0=C_0'=C$. Let $\gamma_i,\gamma_i'$ be the axes of $\rho (C_i),\rho (C_i')$. Let $\eta_i$ and $\eta_i'$ be common orthogonals to $\gamma_{i-1},\gamma_{i+1}$ and $\gamma_{i-1}',\gamma_{i+1}'$, respectively.
The {\it twist-bend parameter} $s(C)$ is the complex distance between unit tangent vector to $\eta_{1}$ at the point $\eta_{1}\cap\gamma_0$ and the unit tangent vector to $\eta'_{1}$ at the point $\eta_{1}'\cap\gamma_0$.
The choices involved guarantee that the complex twist-bend parameter $s(C)$ is well defined in $\mathbb{C}/(2\pi i\mathbb{Z}+hl(C)\mathbb{Z})$.

\section{Pleated surfaces and transverse cocycles to geodesic laminations}

Recall that $S$ is a closed surface of genus $g\geq 2$. Let $\lambda$ be a {\it maximal} geodesic lamination on $S$, namely each component of the complement of $\lambda$ is an ideal hyperbolic triangle.
We do not need to specify a hyperbolic metric on $S$ in order to be able to talk about geodesic laminations on $S$ (cf. \cite{Bon2}). 

Let $\pi :\mathbb{H}^2\to S$ be the universal covering for a
metric $m$ and let $\tilde{\lambda}=\pi^{-1}(\lambda )$. 
An {\it abstract pleated surface} for $S$ with
the pleating locus $\lambda$ is a pleating map $\tilde{f}$ with the
pleating locus $\tilde{\lambda}$ from the hyperbolic plane
$\mathbb{H}^2$ into the hyperbolic three-space $\mathbb{H}^3$ which
is equivariant under the action on $\mathbb{H}^2$ of the covering
group $G$ of $S$ and the action on $\mathbb{H}^3$ of a subgroup
$G_{\tilde{f}}$ of $PSL_2(\mathbb{C})$. If the continuous extension
of $\tilde{f}$ from the ideal boundary
$\partial_{\infty}\mathbb{H}^2$ of $\H^2$ to the ideal boundary
$\partial_{\infty}\mathbb{H}^3$ is injective then the group
$G_{\tilde{f}}$ is quasifuchsian and $\tilde{f}$ projects to a
pleated map from $S=\mathbb{H}^2/G$ into the quasifuchsian
three-manifold $\mathbb{H}^3/G_{\tilde{f}}$.

In this article we consider only finite maximal geodesic laminations on $S$ which are necessarily obtained by triangulating pairs of pants of a pants  decomposition $\mathcal{P}$ of $S$ as follows. Let $\Pi_1$ and $\Pi_2$ be two pairs of pants of
$\mathcal{P}$ that have $C\in\mathcal{P}$ as one of its boundaries.
It is possible that $\Pi_1=\Pi_2$. Assume that ideal triangulations of $\Pi_1$ and $\Pi_2$ are given. Let $a_1^j$ for $j=1,2$ be the
boundary edges of the triangulation of $\Pi_1$ whose one end
accumulate at $C$, and similarly let $a_2^j$ for $j=1,2$ be the
boundary edges of the triangulation of $\Pi_2$ whose one end
accumulate at $C$. Let $\tilde{C}$ be a lift of $C$ to
$\mathbb{H}^2$. Then there are adjacent lifts $\tilde{a}_1^j$ of $a_1^j$ for
$j=1,2$ that share an ideal endpoint $x_1$ with $\tilde{C}$, and
there are adjacent lifts $\tilde{a}_2^j$ of $a_2^j$ for $j=1,2$ that share an
ideal endpoint $x_2$ with $\tilde{C}$. Either $x_1=x_2$ or $x_1\neq
x_2$. If $x_1=x_2$ then, we say that the triangulations of the two pairs of pants
with common boundary $C$ {\it accumulate in the same direction on} $C$.
From now on, we assume that $\lambda_{\mathcal{P}}$ is a finite, maximal geodesic lamination that is obtained by triangulating pairs of pants of $\mathcal{P}$ such that the triangulations of pairs of pants with common boundaries accumulate in the same direction at each $C\in \mathcal{P}$.

 Let $\{\Pi_j\}_{j=1}^{2g-2}$ be the pairs of pants in $\mathcal{P}$. Given  a pair of pants $\Pi_j$ in $\mathcal{P}$ with cuffs $C_i^j\in\mathcal{P}$, $i=1,2,3$,  denote by $\gamma_i^j\in\pi_1(S)$ the elements representing closed curves $C_i^j$ such that $\gamma^j_3\gamma^j_2\gamma^j_1=id$ in $\pi_1(S)$. Let $\{ (hl(C),s(C))\}_{C\in\mathcal{P}}$ be the reduced complex Fenchel-Nielsen coordinates that satisfy (\ref{eq:twist-bend}) and
(\ref{eq:half-length-modif}). By \cite[Proposition 2.3]{Kou}, 
there exists a representation $\rho :\pi_1(S)\to PSL_2(\mathbb{C})$ which realizes  $\{ (hl(C),s(C))\}_{C\in\mathcal{P}}$ such that $\rho (\gamma_i^j\in PSL_2(\mathbb{C})$, for $i=1,2,3$, and $j=1,2,\ldots ,2g-2$ are loxodromic and $\rho (\gamma_i^j)$, $i=1,2,3$, have distinct endpoints for each $j=1,2,\ldots ,2g-2$. 
 Let $S$ be endowed with the hyperbolic metric whose Fenchel-Nielsen coordinates $\{ (Re(hl(C)),Re(s(C)))\}_{C\in\mathcal{P}}$ are the real parts of the reduced complex Fenchel-Nielsen coordinates $\{ (hl(C),s(C)\}_{C\in \mathcal{P}}$.
Let $G$ be the covering group for the universal covering $\pi
:\mathbb{H}^2\to S$. Consider the lifts $\tilde{C}$ to $\mathbb{H}^2$ of the cuffs of $\mathcal{P}$. Then there exists a developing map $\tilde{f}$ from the set of endpoints of $\tilde{C}\in\mathcal{P}$ into $\partial_{\infty}\mathbb{H}^3$ which realizes the reduced complex Fenchel-Nielsen coordinates $\{ (hl(C),s(C)\}_{C\in \mathcal{P}}$. The map $\tilde{f}$ extends to an abstract pleating map as follows. Let $\tilde{\lambda}_{\mathcal{P}}$ be the lift of $\lambda_{\mathcal{P}}$ to $\mathbb{H}^2$. Each geodesic of $\tilde{\lambda}_{\mathcal{P}}$ which is not a lift of a cuff has its both endpoints at the endpoints of two lifts of two different cuffs of a single pair of pants in $\mathcal{P}$ which implies that the endpoints are distinct. Thus $\tilde{f}$ extends to map each geodesic of $\tilde{\lambda}_{\mathcal{P}}$ into a geodesic of $\mathbb{H}^3$. Since the complementary components to $\tilde{\lambda}_{\mathcal{P}}$ are ideal hyperbolic triangles, it follows that we have an extension $\tilde{f}:\mathbb{H}^2\to\mathbb{H}^3$ which determines an abstract pleated  surface with the pleating locus $\tilde{\lambda}_{\mathcal{P}}$. Thus the representation $\rho$ realizes the geodesic lamination $\tilde{\lambda}_{\mathcal{P}}$ (cf. \cite{Bon1}). 

An abstract pleated surface $\tilde{f}:\mathbb{H}^2\to\mathbb{H}^3$ with a pleating locus $\tilde{\lambda}_{\mathcal{P}}$ determines a $(\mathbb{C}/2\pi i\mathbb{Z})$-valued transverse cocycle $\alpha$
to the geodesic lamination $\tilde{\lambda}_{\mathcal{P}}$ (cf. \cite{Bon2}). Namely, $\alpha$ determines a finitely additive assignment of a number in $\mathbb{C}/2\pi i\mathbb{Z}$ to each arc transverse to $\tilde{\lambda}_{\mathcal{P}}$ (with endpoints in the complementary triangles of $\tilde{\lambda}_{\mathcal{P}}$) which is homotopy invariant relative $\lambda_{\mathcal{P}}$. If $k$ is a geodesic arc connecting triangles $\Delta_1$ and $\Delta_2$, then we write $\alpha (\Delta_1,\Delta_2)=\alpha (k)$ because $\alpha (k)$ depends only on the homotopy class of $k$ relative $\tilde{\lambda}_{\mathcal{P}}$.  The real part of $\alpha$ is an $\mathbb{R}$-valued transverse cocycle
which completely determines the path metric on the pleated surface (cf. \cite{Bon2}). The imaginary part of $\alpha$ is an $(\mathbb{R}/2\pi  \mathbb{Z})$-valued transverse cocycle $\beta$ to the geodesic lamination $\lambda_{\mathcal{P}}$. The transverse cocycle $\beta$ determines the amount of the bending of the pleated surface $\tilde{f}:\mathbb{H}^2\to\mathbb{H}^3$ (cf. \cite{Bon2}).

Our first task is to translate the conditions (\ref{eq:half-length-modif}) and (\ref{eq:twist-bend}) in terms of the associated transverse cocycle to $\tilde{\lambda}_{\mathcal{P}}$. 
Let $\tilde{f}:\mathbb{H}^2\to\mathbb{H}^3$ be
the pleating map corresponding to the reduced Fenchel-Nielsen
coordinates $\{ hl(C),s(C)\}_{C\in\mathcal{P}}$ starting from the
real Fenchel-Nielsen coordinates $\{
Re(hl(C)),Re(s(C))\}_{C\in\mathcal{P}}$ as above.  We also lift the decomposition of the pairs of pants of
$\mathcal{P}$ into right-angled hexagons. A right-angled hexagons on
$S$ lifts to an infinite collection of right-angled hexagons in
$\mathbb{H}^2$. Fix a lift $\tilde{C}\in \pi^{-1}(C)$ of the closed
geodesic $C$ and fix a lifted hexagon $\Sigma$ that has one boundary side on
$\tilde{C}$. Then there is a unique lifted hexagon $\Sigma'$ with
one boundary side on $\tilde{C}$ which lies on the opposite side of $\tilde{C}$  
such that the distance between the corresponding vertices
of $\Sigma$ and $\Sigma'$ on $\tilde{C}$ is equal to
$Re(s(C))$. We say that $\Sigma$ and $\Sigma'$ are {\it
$0$-neighbors}. Thus $0$-neighbors hexagons meet along $\tilde{C}$
and have an arc of length $R/2-Re(s(C))$ in common (cf.
Figure 1). If two hexagons $\Sigma$ and $\Sigma''$ meet along their
boundaries but they are not $0$-neighbors then we call them {\it
$1$-neighbors} (cf. Figure 1).

\begin{figure}
\centering
\includegraphics[scale=0.5]{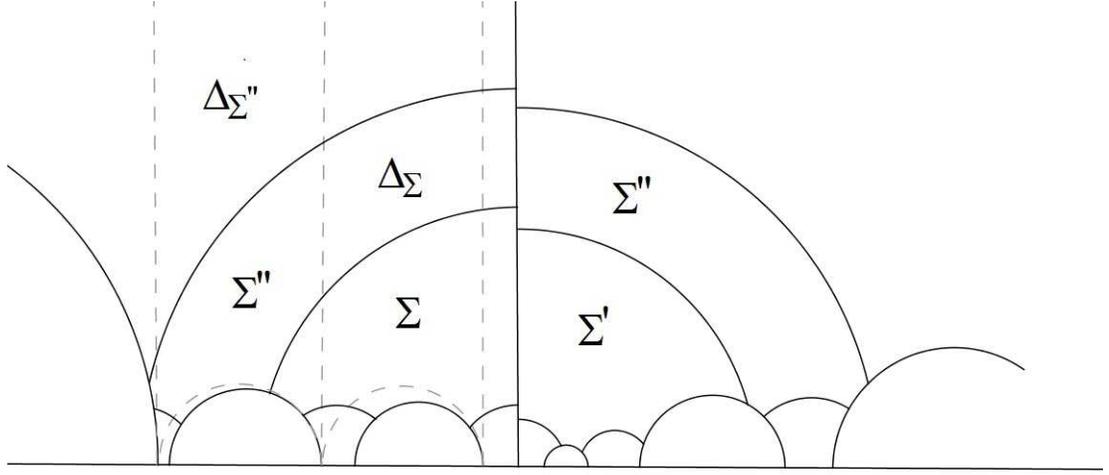}
\caption{$\Sigma$ and $\Sigma'$ are $0$-neighbors; $\Sigma$ and $\Sigma''$ are $1$-neighbors;
$\Delta_{\Sigma}$ is the characteristic triangle of $\Sigma$}
\end{figure}

Fix a lifted hexagon $\Sigma$ in $\mathbb{H}^2$. Among all
complementary triangles to $\tilde{\lambda}_{\mathcal{P}}$ there is
a unique triangle $\Delta_{\Sigma}$ whose all three boundary sides
intersect $\Sigma$. We call $\Delta_{\Sigma}$ the {\it canonical
triangle} for $\Sigma$. Let $\Sigma_t$ be the intersection of
$\Sigma$ and $\Delta_{\Sigma}$. Then $\Sigma\setminus\Sigma_t$ has
three connected components each being a quadrilateral (cf. Figure
1). Let $\mathcal{H}$ be the set of all lifted hexagons in
$\mathbb{H}^2$. Then
$$
\mathcal{TH}_t=\bigcup_{\Sigma\in\mathcal{H}}\Sigma_t
$$
separates geodesics in $\pi^{-1}(\mathcal{C})$ (cf. Figure 2). We give a proof of Theorem \ref{thm:f-n into shear-bend} from Introduction.

\begin{figure}
\centering
\includegraphics[scale=0.5]{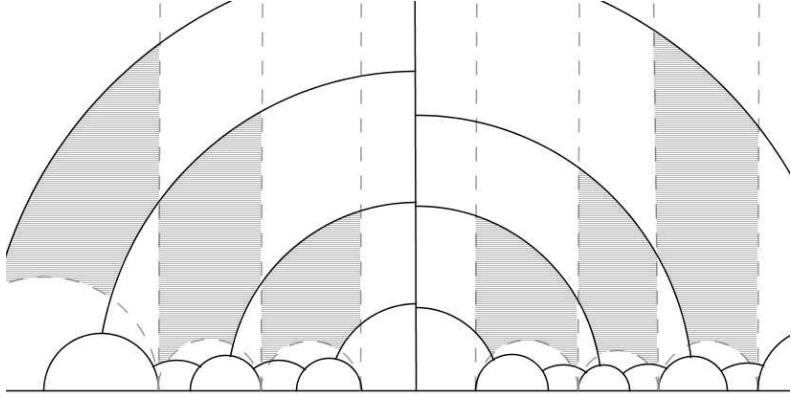}
\caption{The set $\mathcal{TH}_t$}
\end{figure}

\begin{theorem}
\label{lem:f-n into shear-bend} There exists $C_0>0$ such that the following holds. Let $\{
hl(C),s(C)\}_{C\in\mathcal{P}}$ be the reduced complex Fenchel-Nielsen
coordinates such that
$$
|hl(C)-R/2|<\frac{\epsilon}{R}
$$
and 
$$
|s(C)-1|<\frac{\epsilon}{R}.
$$
Let $S$ be endowed with a hyperbolic metric whose
Fenchel-Nielsen coordinates are $\{
Re(hl(C)),Re(s(C))\}_{C\in\mathcal{P}}$. Let
$\lambda_{\mathcal{P}}$ be a maximal geodesic lamination obtained by
triangulating pairs of pants of $\mathcal{P}$ such that the edges of
the triangles from both sides of each $C\in\mathcal{P}$ accumulate
in the same direction. Let
$\tilde{f}:\mathbb{H}^2\to\mathbb{H}^3$ be the bending map with the
bending locus
$\tilde{\lambda}_{\mathcal{P}}=\pi^{-1}(\lambda_{\mathcal{P}})$
which realizes the complex Fenchel-Nielsen coordinates $\{
hl(C),s(C)\}_{C\in\mathcal{P}}$, where $\pi :\mathbb{H}^2\to S$ is
the universal covering. Denote by $\beta$ the bending cocycle
transverse to $\tilde{\lambda}_{\mathcal{P}}$ for $\tilde{f}$. Let
$\tilde{l}$ be an isolated leaf of $\tilde{\lambda}_{\mathcal{P}}$ which is
on the boundary of two complementary triangles $\Delta_1(\tilde{l})$ and $\Delta_2(\tilde{l})$ of
$\tilde{\lambda}_{\mathcal{P}}$. Then
$$
|\beta (\Delta_1(\tilde{l}),\Delta_2(\tilde{l}))|\leq \frac{C_0\epsilon}{R} .
$$
Moreover, let $\Sigma_1$ and $\Sigma_2$ be $0$-neighbor hexagons and
let $\Delta_{\Sigma_1}$ and $\Delta_{\Sigma_2}$ be their canonical
triangles, respectively. Then
$$
|\beta (\Delta_{\Sigma_1},\Delta_{\Sigma_2})|\leq
\frac{C_0\epsilon}{R}.
$$
\end{theorem}

\begin{proof}
Let $\Pi$ be a pair of pants in $\mathcal{P}$ with boundary curves
$C_i$, $i=0,1,2$. Let $l_i$, for $i=0,1,2$, be the geodesics of $\lambda_{\mathcal{P}}$ that triangulate
$\Pi$ such that $l_{i+1}$ and $l_{i+2}$ accumulate on $C_i$, for $i=0,1,2$, with the indices taken modulo $3$. Let $\tilde{C}_i$ be a lift of $C_i$ to $\mathbb{H}^2$ and let $\tilde{l}_{i+1},\tilde{l}_{i+2}$ be consecutive lifts of $l_{i+1},l_{i+2}$ that share a common endpoint with $\tilde{C}_i$. Let $\gamma_i\in PSL_2(\mathbb{R})$ be the deck transformation corresponding to $\tilde{C}_i$ such that $\gamma_i(\tilde{l}_{i+1})=\tilde{l}'_{i+1}$ is adjacent to $\tilde{l}_{i+2}$. Let $\tilde{l}_{i+2}'=\gamma_i(\tilde{l}_{i+2})$. Then $\tilde{l}'_{i+2}$ is adjacent to $\tilde{l}_{i+1}'$.

Let $r_{i+2}$ be the geodesic ray orthogonal to $\tilde{f}(\tilde{l}_{i+2})$ that starts at the endpoint of $\tilde{f}(\tilde{l}_{i+1})$ which is not in common with $\tilde{f}(\tilde{C}_i)$. Let $r_{i+2}'$ be the geodesic ray orthogonal to $\tilde{f}(\tilde{l}_{i+2})$ that starts at the endpoint of $\tilde{f}(\tilde{l}_{i+1}')$ which is not in common with $\tilde{f}(\tilde{C}_i)$. Define $s_{i+2}$ to be the complex distance between the unit tangent vector to $r_{i+2}'$ at $r_{i+2}'\cap\tilde{f}(\tilde{l}_{i+2})$ and the unit tangent vector to $r_{i+2}$ at $r_{i+2}\cap\tilde{f}(\tilde{l}_{i+2})$.  Define $s_{i+1}$ using the geodesics $\tilde{f}(\tilde{l}_{i+2}),\tilde{f}(\tilde{l}_{i+1}'),\tilde{f}(\tilde{l}_{i+2})$ similar to the above. Then
$$
s_{i+1}+s_{i+2}=l(\delta_i)
$$
where $\delta_i=\tilde{f}\circ\gamma_i\circ\tilde{f}^{-1}\in PSL_2(\mathbb{C})$ and $l(\delta_i)$ is the complex translation length of $\delta_i$, for $i=0,1,2$. Solving the above system gives
$$
s_i=\frac{l(\delta_{i+1})+l(\delta_{i+2})-l(\delta_i)}{2}
$$
for $i=0,1,2$. Since 
$$\big{|}Im(\frac{1}{2}l(\delta_i))\big{|}\leq\frac{\epsilon}{R},$$ it
follows that
$$
|\beta (\Delta_1(\tilde{l}_i),\Delta_2(\tilde{l}_i))|\leq \frac{C_0\epsilon}{R}
$$
for $i=0,1,2$ and some $C_0>0$.

Let $\Sigma\subset\mathbb{H}^2$ be a lifted right angled hexagon from a pair of
pants $\Pi$ whose boundary sides lie on $C_i$, $i=0,1,2$. We fix lifts
$\tilde{C}_i$ of $C_i$ such that three sides of
$\Sigma$ lie on $\tilde{C}_i$ for $i=0,1,2$. Let $\Sigma_{\tilde{f}}$ be the skew right angled
hexagon whose three sides lie on $\tilde{f}(C_1)$, $\tilde{f}(C_2)$
and $\tilde{f}(C_3)$. Then the complex lengths of these sides are
$hl(C_1)$, $hl(C_2)$ and $hl(C_3)$. These sides are called {\it long
sides} and the other three sides of $\Sigma_{\tilde{f}}$ are called
{\it short sides}. Denote by $h_i$ the short side of $\Sigma_{\tilde{f}}$
which connects $\tilde{f}(C_{i+1})$ and $\tilde{f}(C_{i+2})$. Then
the hexagon cosine formula directly gives (cf. \cite{Bow},
\cite{KahnMark}))
$$
l(h_i)=2e^{-\frac{R}{4}+\frac{1}{2}[hl(C_{i+1})+hl(C_{i+2})-hl(C_i))]}+O(e^{-3R/4})
$$
where $l(h_i)$ is the complex distance between
$\tilde{f}(\tilde{C}_{i+1})$ and $\tilde{f}(\tilde{C}_{i+2})$. This implies
$$
|Re(l(h_i))|=O\big{(}e^{-R/4}\big{)}
$$
and
$$
|Im(l(h_i))|=O\big{(}\frac{\epsilon}{R}e^{-R/4}\big{)}.
$$

Let $\Sigma_1$ and $\Sigma_2$ be two $0$-neighbors hexagons in
$\mathbb{H}^2$, and let $\Delta_{\Sigma_1}$ and $\Delta_{\Sigma_2}$
be their canonical triangles. Let
$\tilde{C}\in\tilde{\mathcal{P}}=\pi^{-1}(\mathcal{P})$ be the
geodesic which separates $\Delta_{\Sigma_1}$ and
$\Delta_{\Sigma_2}$, and let $C=\pi (\tilde{C})\in\mathcal{P}$. Note that both $\Sigma_1$ and $\Sigma_2$ have one boundary side on $\tilde{C}$.
Normalize the bending map such that the ideal triangles
$\tilde{f}(\Delta_{\Sigma_1})$ and $\tilde{f}(\Delta_{\Sigma_2})$,
have a common endpoint $\infty$, and that $\tilde{f}(\tilde{C})$ has
endpoints $0$ and $\infty$. Let $\tilde{C}_1^j$ for $j=1,2$ be the two
geodesics of $\tilde{\mathcal{P}}$ (different from $\tilde{C}$) which contain boundary sides of
$\Sigma_1$, and let $\tilde{C}_2^j$ for $j=1,2$ be the two geodesics
in $\tilde{\mathcal{P}}$ (different from $\tilde{C}$) which contain boundary sides of $\Sigma_2$.
We can assume that the twist-bend $s(C)$ is the complex distance (along $\tilde{f}(\tilde{C})$)
between the common orthogonal to $\tilde{f}(\tilde{C}_1^1)$ and $\tilde{f}(\tilde{C})$, and the
common orthogonal to $\tilde{f}(\tilde{C}_2^1)$ and $\tilde{f}(\tilde{C})$. It follows that the
complex distance (along $\tilde{f}(\tilde{C})$) between the common orthogonal to $\tilde{f}(\tilde{C}_1^2)$ and
$\tilde{f}(\tilde{C})$, and the common orthogonal to $\tilde{f}(\tilde{C}_2^2)$ and $\tilde{f}(\tilde{C})$ is
also equal to the twist-bend $s(C)$ (cf. Figure 3).

\begin{figure}
\centering
\includegraphics[scale=0.5]{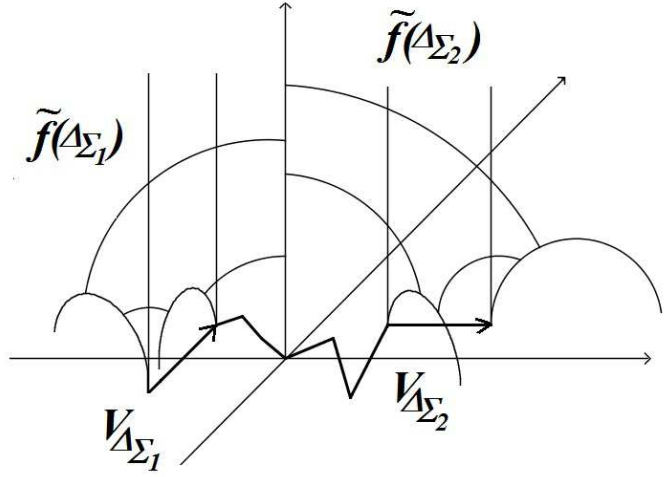}
\caption{Computing the bending cocycle.}
\end{figure}

We recall the definition of $\beta
(\Delta_{\Sigma_1},\Delta_{\Sigma_2})$ given by Bonahon \cite{Bon2}.
Let $\mathcal{W}$ be the component of $\mathbb{H}^2\setminus
(\Delta_{\Sigma_1}\cup\Delta_{\Sigma_2})$ which separates
$\Delta_{\Sigma_1}$ and $\Delta_{\Sigma_2}$. Denote by
$\tilde{\lambda}_{\mathcal{P}}(\Delta_{\Sigma_1},\Delta_{\Sigma_2})$
the set of leaves of $\tilde{\lambda}_{\mathcal{P}}$ which separate
$\Delta_{\Sigma_1}$ and $\Delta_{\Sigma_2}$, and orient them to the
left as seen from $\Delta_{\Sigma_1}$. The leaves of
$\tilde{\lambda}_{\mathcal{P}}$ divide $\mathcal{W}$ into hyperbolic
strips and the images under $\tilde{f}$ of the hyperbolic strips are
two-dimensional hyperbolic strips in $\mathbb{H}^3$. Each such
hyperbolic strip intersects $\partial_{\infty}\mathbb{H}^3$ in two
circular arcs with a possibility that one is reduced to a point such
that one circular arc is bounded by the negative endpoints of the
leaves $\tilde{f}(\tilde{\lambda}_{\mathcal{P}}(\Delta_{\Sigma_1},
\Delta_{\Sigma_2}))$, and the other by positive endpoints. Let $\gamma\in\partial_{\infty}\mathbb{H}^3$ be an oriented,
piecewise circular curve 
formed by concatenating the circular arcs bounded by negative endpoints from $\Delta_{\Sigma_1}$ to $\Delta_{\Sigma_2}$. Let
$v_{\Delta_{\Sigma_1}}$ be the outward tangent vector to the circular arc of
the intersection of
$\tilde{f}(\Delta_{\Sigma_1})\cap\partial_{\infty}\mathbb{H}^3$, and let $v_{\Delta_{\Sigma_2}}$ be the inward tangent
vector to the circular arc of the intersection
$\tilde{f}(\Delta_{\Sigma_2})\cap\partial_{\infty}\mathbb{H}^3$. Then (cf. \cite{Bon2})
$$
\beta (\Delta_{\Sigma_1},\Delta_{\Sigma_2})=\angle
(v_{\Delta_{\Sigma_1}},v_{\Delta_{\Sigma_2}})-\sum_{W}\beta_W
$$
where $\angle
(v_{\Delta_{\Sigma_1}},v_{\Delta_{\Sigma_2}})$ is the angle under Euclidean parallel transport in $\mathbb{C}$,  $\beta_W$ is the signed curvature of the circular subarc $W$
of $\gamma$ and the sum is over all circular subarcs of $\gamma$. In
our case, all circular arcs are Euclidean segments and each term of the sum in
the above formula is zero. Thus we obtain
$$
\beta (\Delta_{\Sigma_1},\Delta_{\Sigma_2})=\angle
(v_{\Delta_{\Sigma_1}},v_{\Delta_{\Sigma_2}}).
$$

To finish the proof we refer to Figure 3. The vector
$v_{\Delta_{\Sigma_1}}$ is parallel to the vector
$\overrightarrow{xy}$ in $\partial_{\infty}\mathbb{H}^3$ whose
initial point $x$ is an endpoint of $\tilde{f}(\tilde{C}^1_1)$ and terminal point $y$
is an endpoint of $\tilde{f}(\tilde{C}^1_2)$ as in Figure 3. Similarly, the vector
$v_{\Delta_{\Sigma_2}}$ is parallel to the vector
$\overrightarrow{y'x'}$ in $\partial_{\infty}\mathbb{H}^3$ whose
initial point $y'$ is an endpoint of $\tilde{f}(\tilde{C}^2_2)$ and terminal point $x'$
is an endpoint of $\tilde{f}(\tilde{C}^2_1)$ as in Figure 3. We normalize the situation
such that the short sides $h_1$ and $h_2$ of $\Sigma_1$ meet $\tilde{C}$ at
$j=(0,0,1)\in\mathbb{H}^3$ and $e^{-Re(hl(C))}j=Ce^{-R/2}j\in\mathbb{H}^3$
for $e^{-\frac{\epsilon}{R}}<C< e^{\frac{\epsilon}{R}}$. In this case, the points where
the short sides $h_1'$ and $h_2'$ of $\Sigma_2$ meet $\tilde{C}$ are
$e^{Re(s(C))}j=Cej\in\mathbb{H}^3$,
$e^{-\frac{\epsilon}{R}}<C<e^{\frac{\epsilon}{R}}$, and
$e^{Re(s(C)-hl(C))}j=C_1e^{1-\frac{R}{2}}j\in \mathbb{H}^3$,
$e^{-\frac{2\epsilon}{R}}<C_1<e^{\frac{2\epsilon}{R}}$.

If $h_1$ lies in the $xz$-plane in $\mathbb{H}^3$ then $x$ is an analytic function of the complex length $l(h_1)$ of $h_1$. An explicit (and elementary) computation shows that the derivative of $x$ in the variable $l(h_1)$ at the point $l(h_1)=0$ is non-zero. Thus the  euclidean distance from $x$ to $0\in\partial_{\infty}\mathbb{H}^3$ is $O(|l(h_1)|)=O(e^{-R/4})$ (this holds without the restriction that $h_1$ is in the $xz$-plane). Since $y$ is the image of $x$ under the map $z\mapsto e^{-hl(C)}z$, it follows that the distance between $y$ and $0$ is $O(e^{-3R/4})$. Similar statements hold for $x'$ and $y'$, respectively. Consider the Fenchel-Nielsen coordinates $\{ (Re(hl(C)),s(C)\}_{C\in\mathcal{P}}$ and let $\tilde{f}_{Re}$ be the corresponding developing map. We normalize $\tilde{f}_{Re}$ such that $\tilde{f}_{Re}(\tilde{C})$ has endpoints $0$ and $\infty$, and that $\tilde{f}_{Re}(\Delta_{\Sigma_1})$ and $\tilde{f}_{Re}(\Delta_{\Sigma_2})$ have a common endpoint $\infty$. Moreover, we require that the common orthogonal between $\tilde{f}_{Re}(\tilde{C})$ and $\tilde{f}_{Re}(\tilde{C}_1^1)$ meets $\tilde{f}_{Re}(\tilde{C})$ at $j\in\mathbb{H}^3$, and the common orthogonal between  $\tilde{f}_{Re}(\tilde{C})$ and $\tilde{f}_{Re}(\tilde{C}_1^2)$ meets $\tilde{f}_{Re}(\tilde{C})$ at $e^{1+O(\epsilon /R)}j\in\mathbb{H}^3$. Let $x_0$, $y_0$, $x_0'$ and $y_0'$ be the endpoints of $\tilde{f}_{Re}(\tilde{C}_1^1)$, $\tilde{f}_{Re}(\tilde{C}_2^1)$, $\tilde{f}_{Re}(\tilde{C}_1^2)$ and $\tilde{f}_{Re}(\tilde{C}_2^2)$ that are different from $\infty$, respectively.

Let $\tilde{f}_{Im}$ be the developing map which maps the pleated surface for $\{ (Re(hl(C)),s(C))\}_{C\in\mathcal{P}}$ to the pleated surface for $\{ (hl(C),s(C))\}_{C\in\mathcal{P}}$ and fixes $\tilde{C}$. Then $x=\tilde{f}_{Im}(x_0)$,  $y=\tilde{f}_{Im}(y_0)$,  $x'=\tilde{f}_{Im}(x_0')$ and  $y'=\tilde{f}_{Im}(y_0')$. In terms of the geometry, $x$ is the image of the endpoint $x_0$ of the geodesic $\tilde{f}_{Re}(\tilde{C}_1^1)$ under the rotation around the common orthogonal to $\tilde{f}_{Re}(\tilde{C})=\tilde{C}$ and $\tilde{f}_{Re}(\tilde{C}_1^1)$ with the angle of the rotation equal to the imaginary part of the complex length of the common orthogonal to $\tilde{f}(\tilde{C})=\tilde{C}$ and $\tilde{f}(\tilde{C}_1^1)$. The cosine formula estimates this angle to be $O(\frac{\epsilon}{R}e^{-R/4})$. Then the euclidean distance between $x_0$ and $x$ is $O(\frac{\epsilon}{R}e^{-R/4})$ which implies that the distance between $x$ and$0$ is $O(e^{-R/4})$. Note that $y_0$ and $y$ are the images of $x_0$ and $x$ under the maps $z\mapsto e^{-R/4+(1+i)O(\frac{\epsilon}{R})}$. Thus the euclidean distance between $y$ and $y_0$ is $O(\frac{\epsilon}{R}e^{-3R/4})$ and the distance between $y$ and $0$ is $O(e^{-3R/4})$. Similar properties hold for $x_0'$, $x'$ and $y_0'$, $y'$. The angle between the vectors $\overrightarrow{x_0y_0}$ and $\overrightarrow{y_0'x_0'}$ is $s(C)$ because the length of $C\in\mathcal{P}$ is real (which means that the hexagons are not skewed).
The above shows that the angle $\overrightarrow{x_0y_0}$ and $\overrightarrow{xy}$ is $O(\frac{\epsilon}{R})$, and the same estimate for the angle between  $\overrightarrow{y_0'x_0'}$ and $\overrightarrow{y'x'}$. Thus the angle between $\overrightarrow{xy}$ and $\overrightarrow{y'x'}$ is $O(\frac{\epsilon}{R})$ which finishes the proof.
\end{proof}

\subsubsection*{The definition of the bending map}
Let $\Delta_1$ and $\Delta_2$ be two complementary triangles to
$\tilde{\lambda}$. Bonahon \cite{Bon2} defined the bending map $\tilde{f}_{\beta}|_{\Delta_2}=\varphi_{\Delta_1,\Delta_2}$ normalized to be the identity at $\Delta_1$ as follows.  
Let $\mathcal{P}_p=\{\Delta_1',\Delta_2',\ldots ,\Delta_p'\}$ be a sequence of complementary triangles to $\tilde{\lambda}$ which separate $\Delta_1$ and $\Delta_2$ given in the order from $\Delta_1$ to $\Delta_2$. Define
$$
\psi_p=R_{g_{\Delta_1'}^{\Delta_1}}^{\beta(\Delta_1,\Delta_1')}\circ 
R_{g_{\Delta_1'}^{\Delta_2}}^{-\beta(\Delta_1,\Delta_1')}\circ R_{g_{\Delta_2'}^{\Delta_1}}^{\beta(\Delta_1,\Delta_2')}\circ 
R_{g_{\Delta_2'}^{\Delta_2}}^{-\beta(\Delta_1,\Delta_2')}\circ \cdots \circ
R_{g_{\Delta_p'}^{\Delta_1}}^{\beta(\Delta_1,\Delta_p')}\circ 
R_{g_{\Delta_p'}^{\Delta_2}}^{-\beta(\Delta_1,\Delta_p')}
$$
where $R_g^b$ is the hyperbolic rotation around the axis $g\subset\mathbb{H}^3$ by the angle $b\in\mathbb{R}$, and $g_{\Delta_i'}^{\Delta_k}$ is the geodesic on the boundary of $\Delta_i'$ which is closest to $\Delta_k$ for $k=1,2$.  Let $\mathcal{P}$ be the family of all complementary triangles to $\tilde{\lambda}$ that separate $\Delta_1$ and $\Delta_2$. If $\mathcal{P}_p\to\mathcal{P}$ in the sense that $\mathcal{P}_p$ is an increasing family with $\cup_{p=1}^{\infty}\mathcal{P}_p=\mathcal{P}$, then the limit
$$
\psi_{\Delta_1,\Delta_2}=\lim_{\mathcal{P}_p\to\mathcal{P}}\psi_p
$$
exists and it is independent of the choice of $\mathcal{P}_p$ (cf. \cite{Bon2}). Then
$$
\varphi_{\Delta_1,\Delta_2}=\psi_{\Delta_1,\Delta_2}\circ R_{g_{\Delta_2}^{\Delta_1}}^{\beta (\Delta_1,\Delta_2)}.
$$

\vskip .2 cm

The following lemma is established in \cite{KahnMark}. We give a
different proof below.

\begin{lemma}
\label{lem:length1-intersects-at-most-R} Under the above assumptions, a geodesic arc in
$\mathbb{H}^2$ of length $1$ intersects at most $2R+2$ geodesics
from $\pi^{-1}(\mathcal{P})$, when $R$ is large enough.
\end{lemma}

\begin{proof}
Let $l$ be an arc of length $1$ which transversely intersects
geodesics of $\tilde{\mathcal{P}}=\pi^{-1}(\mathcal{P})$.
Let $\{ \tilde{C}_1,\tilde{C}_2,\ldots ,\tilde{C}_n\}$
be the geodesics in $\tilde{\mathcal{P}}$ which intersect $l$ in the given
order and we orient them to the left as seen from the half-plane in
$\mathbb{H}^2\setminus\tilde{C}_1$ which does not contain
$\tilde{C}_2$.
For $R$ large enough, consecutive geodesics in $\{ \tilde{C}_1,\tilde{C}_2,\ldots ,\tilde{C}_n\}$ are connected by the short arcs of the hexagons (otherwise $l$ would intersect two short sides of a single hexagon which would imply $|l|\geq R/4>1$).
 Given $\tilde{C}_j$ and $\tilde{C}_{j+1}$, let $h_j$
be the common orthogonal, and let $x_j^{+}=\tilde{C}_j\cap h_j$ and
$x_{j+1}^{-}=\tilde{C}_{j+1}\cap h_j$ for $j=1,\ldots ,n-1$, and
$x_n^{+}=x_n^{-}$. Given $a\in \tilde{C}_j$, define $r(a)$ to be the
signed distance between $a$ and $x^{+}_j$. Issue a geodesic $g_a$
through $a$ such that the angle of intersection between
$\tilde{C}_j$ and $g_a$ is equal to the angle of intersection
between $\tilde{C}_{j+1}$ and $g_a$. Let
$a'=g_a\cap\tilde{C}_{j+1}$. Then the signed distance between $a'$
and $x_{j+1}^{-}$ is equal to $r(a)$ and consequently the signed
distance between $a'$ and $x_{j+1}^{+}$ is $r(a)-(1\pm\frac{\epsilon}{R})$ because the twist parameter is
$1\pm\frac{\epsilon}{R}$ by the assumption. Let $r^{-}(a)$ denote the signed distance
of $a$ to $x_j^{-}$ for $a\in C_j$. Thus $r^{-}(a')=r(a)$ for $a\in
C_j$.

\begin{figure}
\centering
\includegraphics[scale=0.5]{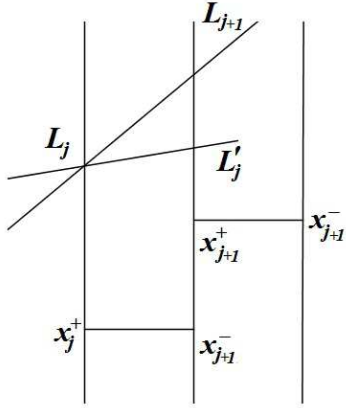}
\caption{The number of intersections.}
\end{figure}

Let $L_j=l\cap\tilde{C}_j$. We compare the signed distance between
$L_{j+1}$ and $x_{j+1}^{+}$ to the signed distance between $L_j'$
and $x_{j+1}^{+}$. Consider the hyperbolic triangle with vertices
$L_j$, $L_j'$ and $L_{j+1}$. The angle at $L_j$ is smaller than the
angle at $L_j'$ (cf. Figure 4). By the sine formula for hyperbolic
triangles we get
\begin{equation}
\label{eq:dist_l}
d(L_j,L_{j+1})>d(L_j',L_{j+1}).
\end{equation}
By the condition on the twist parameters we have
$$
r(L_{j+1})=r^{-}(L_{j+1})-(1\pm\frac{\epsilon}{R}).
$$
Moreover, we have
$$
r^{-}(L_{j+1})=r^{-}(L_j')\pm d(L_j',L_{j+1})=r(L_j)\pm
d(L_j',L_{j+1}).
$$
The above gives
$$
r(L_{j+1})=r(L_j)\pm d(L_j',L_{j+1})-(1\pm\epsilon )
$$
and thus
\begin{equation}
\label{eq:dist_l1}
r(L_n)\geq r(L_1)- 1+ (n-1)(1\pm\frac{\epsilon}{R})
\end{equation}
for all $n\geq 2$ because $\sum_{i=1}^n d(L_j',L_{j+1})\leq\sum_{i=1}^n d(L_j,L_{j+1})\leq 1$ by (\ref{eq:dist_l}).

Assume that $n\geq 2R+2$. Then there is $1\leq j\leq n$ such that
$|r(L_j)|\geq R$ by (\ref{eq:dist_l1}). We find the contradiction with this inequality by
proving that $d(L_j,L_{j+1})$ is too large in this case.

We prove that $d(L_j,L_{j+1})$ is too large. Without the loss of generality, we assume that the quadrilateral $Q$ with vertices $L_j$, $x_j^{+}$, $x_{j+1}^{-}$ and $L_{j+1}$ has a right angle at the vertex $L_{j+1}$. It follows then that thee angles of $Q$ are equal to $\frac{\pi}{2}$ since the angles at $x_j^{+}$ and $x_{j+1}^{-}$ are equal to $\frac{\pi}{2}$. An elementary hyperbolic geometry gives
$$
\cosh^2 d(L_j,L_{j+1})=\cosh^2 d(x_j^{+},L_j)\sinh^2 d(x_j^{+},x_{j+1}^{-})+1.
$$
Since 
$$
d(x_j^{+},L_j)=|r(L_j)|\geq R
$$
and 
$$
d(x_j^{+},x_{j+1}^{-})\leq Ce^{-R/2},
$$
the above gives
$$
d(L_j,L_{j+1})\geq\frac{R}{2}-C
$$
for a fixed $C>0$ and $R$ large enough. This implies that $d(L_j,L_{j+1})>1$ for $R$ large enough which is a contradiction. Thus a geodesic arc of length $1$ intersects at most $2R+2$ geodesics of $\tilde{\mathcal{P}}$.
\end{proof}

\section{Injectivity of the bending maps}

The purpose of this section is to prove the following theorem which is the first statement of Theorem \ref{thm:injective-bending-cocycles} from Introduction. We finish the proof of the remaining statements of Theorem \ref{thm:injective-bending-cocycles} in the next section.

\begin{theorem}
\label{tim:half-of-main}
Given $C_0>0$, there exist $\hat{\epsilon}>0$ and $R(\hat{\epsilon})>0$ such that for each $0\leq\epsilon <\hat{\epsilon}$ and $R\geq R(\hat{\epsilon})$ the following is satisfied. Let $S$ be a closed hyperbolic surface equipped with a maximal, finite geodesic lamination $\lambda$ such that each closed geodesic of $\lambda$ has length in the interval $(R-\frac{\epsilon}{R},R+\frac{\epsilon}{R})$ and that each geodesic arc of length $1$ intersects at most $C_0 R$ closed geodesics of $\lambda$.  Assume that a bending cocycle $\beta$ transverse to the lift $\tilde{\lambda}$ in $\mathbb{H}^2$ satisfies
\begin{equation}
\label{eq:isolated}
|\beta (\Delta^1(\tilde{l}),\Delta^2(\tilde{l}))|\leq \frac{C_0\epsilon}{R} 
\end{equation}
for each isolated leaf $\tilde{l}$ and complementary triangles $\Delta^1(\tilde{l})$ and $\Delta^2(\tilde{l})$ with common boundary $\tilde{l}$, and
\begin{equation}
\label{eq:cuff}
|\beta (\Delta_{\Sigma_1},\Delta_{\Sigma_2})|\leq
\frac{C_0\epsilon}{R}
\end{equation}
for the characteristic triangles $\Delta_{\Sigma_1}$ and $\Delta_{\Sigma_2}$ of each two $0$-neighbors hexagons $\Sigma_1$ and $\Sigma_2$ (coming from the pants decomposition of $S$ whose cuffs are closed geodesics of $\lambda$). Then the induced bending map
$$
\tilde{f}_{\beta}:\partial_{\infty}\mathbb{H}^2\to\partial_{\infty}\mathbb{H}^3
$$
is injective. 
\end{theorem}

\begin{proof}
Let $x,y\in\partial_{\infty}\mathbb{H}^2$ be two different points.
We need to prove that $\tilde{f}_{\beta}(x)\neq\tilde{f}_{\beta}(y)$. Let $\mathcal{P}$ be a pants decomposition of $S$ whose cuffs are closed curves of $\lambda$. We fix a decomposition of $S$ into hexagons as in \S 3 using the pants decomposition $\mathcal{P}$ and lift it to the universal covering $\pi :\mathbb{H}^2\to S$. Recall that $\tilde{\lambda}=\pi^{-1}(\lambda )$ and $\tilde{\mathcal{P}}=\pi^{-1}(\mathcal{P})$.  Let $g$ be the
geodesic in $\mathbb{H}^2$ whose ideal endpoints are $x$ and $y$. If
$g\subset \mathbb{H}^2\setminus\mathcal{TH}_t$ then $g$ is a lift of
some $C\in\mathcal{P}$ and $\tilde{f}_{\beta}(x)\neq\tilde{f}_{\beta}(y)$ because
$g$ is in the bending locus of $\tilde{f}_{\beta}$.

Therefore we assume that $g\cap \mathcal{TH}_t\neq\emptyset$. Fix a
hexagon $\Sigma^{0}$ such that $g\cap\Sigma_t^{0}\neq\emptyset$. Let
$P$ be a point in $g\cap\Sigma_t^{0}$. The point $P$ divides the
geodesic $g$ into two rays $g_{\pm 1}$. Let $P_0=P$ and
assume that we have chosen points $P_{\pm 1},P_{\pm 2},\ldots ,P_{\pm n}$ in
the increasing order on $g_{\pm 1}$ such
that $P_{\pm k}\in (\Sigma_{\pm k})_t$ for distinct hexagons $\Sigma_{\pm k}$, for 
$k=1,2,\ldots ,n$. We define $P_{\pm (n+1)}$ as follows.
Let $\Sigma_{\pm (n+1)}$ be the first hexagon after $\Sigma_{\pm n}$ such
that $(\Sigma_{\pm (n+1)})_t$ intersect $g_{\pm 1}$ and that there exists a
point $P_{\pm (n+1)}\in g_{\pm 1}\cap (\Sigma_{\pm (n+1)})_t$ with
$d(P_{\pm n},P_{\pm (n+1)})\geq 1$. If such hexagon does not exist,
then we set $P_{\pm (n+1)}$ to be the ideal endpoint of $g_{\pm 1}$. In this
fashion we partition each $g_{\pm 1}$ into consecutive arcs of lengths at
least $1$. It is possible that the partition is finite when
$P_{\pm (n+1)}$ is the endpoint of $g_{\pm 1}$.

Let $g\subset\mathbb{H}^3$ be a geodesic ray with initial point
$p_0$, and let $p\in g$ be another point. For $0<\theta <\pi$, the {\it
cone $\mathcal{C}(p,g,\theta )$ with vertex $p$, axis $g$ and angle
$\theta$} is the set of all $w\in\mathbb{H}^3$ such that the angle
at $p$ between the positive direction of $g$ and the geodesic ray
from $p$ through $w$ is less than $\theta$. Note that a cone is an open set. A non-zero vector
$(p,v)\in T^1(\mathbb{H}^3)$ uniquely determines a geodesic ray
$g$ which starts at the basepoint $p$ of $v$ and which is tangent to
$v$. Then $\mathcal{C}(p,v,\theta )$ is by the definition
$\mathcal{C}(p,g,\theta )$. The {\it shadow of the cone}
$\mathcal{C}(p,g,\theta )$ is the set
$\partial_{\infty}\mathcal{C}(p,g,\theta )$ of endpoints at
$\partial_{\infty}\mathbb{H}^3$ of all geodesic rays starting at $p$
and inside $\mathcal{C}(p,g,\theta )$. The shadow of a cone is an open subset of $\partial_{\infty}\mathbb{H}^3$.

For $d>0$, let $p_d\in g$ be the point on $g$ which is on the
distance $d$ from $p_0=p$. Let $\eta >0$ be the maximal angle such
that $\mathcal{C}(p_d,g,\eta )\subset\mathcal{C}(p_0,g,\theta )$.
Then $\eta =\eta (d,\theta )$ is a continuous function of $d$ and
$\theta$. For a fixed $0<\theta <\pi$, we have $\eta (d,\theta
)>\theta$ and $\eta (d,\theta )\to\theta$ as $d\to 0$. These
properties are elementary.

Let $\{ P_{\pm n}\}_{n}$ be the points of the partition of $g_{\pm 1}$. We
consider a sequence of cones $\{\mathcal{C}(P_{\pm n},g_{\pm 1},\frac{\pi}{2})\}$. Then 
$$\overline{\partial_{\infty}(\mathcal{C}(P_{\pm (n+1)},g_{\pm 1},\frac{\pi}{2})}\subset \partial_{\infty}
(\mathcal{C}(P_{\pm n},g_{\pm 1},\frac{\pi}{2})
$$
for each $n\in\mathbb{N}$ and we say that the sequence of cones is {\it nested}.

If we prove that
the images of the nested cones under the bending map $\tilde{f}_{\beta}$ remain nested then we are done. Indeed, since $x$ and $y$ lie in the intersection of the shadows of all nested cones along $g_1$ and $g_{-1}$, since the shadows of $\mathcal{C}(P_0,g_1,\frac{\pi}{2})$ and $\mathcal{C}(P_0,g_{-1},\frac{\pi}{2})$ are disjoint, and if $\tilde{f}_{\beta}$ preserves the nesting of the cones, it follows that $\tilde{f}_{\beta}(x)\neq\tilde{f}_{\beta}(y)$. It remains to prove that $\tilde{f}_{\beta}$ preserves the nesting of the cones.  
To see this, it is
enough to normalize $\tilde{f}_{\beta}$ to be the identity on the canonical
triangle $\Delta_{\Sigma_{\pm n}}$ of $\Sigma_{\pm n}$ and to prove that
$$\overline{(\tilde{f}_{\beta}|_{\Delta_{\Sigma_{\pm (n+1)}}})(\partial_{\infty}
\mathcal{C}(P_{\pm (n+1)},g_{\pm 1},\frac{\pi}{2}))}
\subset\partial_{\infty}\mathcal{C}(P_{\pm n},g_{\pm 1},\frac{\pi}{2})$$ for
each $n\in\mathbf{N}$.

Let $a_{\pm n}$ be the arc of $g_{\pm 1}$ between $P_{\pm n}$ and $P_{\pm (n+1)}$.
Note that the length of $a_{\pm n}$ is at least $1$ and that it can be
infinite. We first assume that $a_{\pm n}$ has finite length. 
Let $\Sigma_1,\Sigma_2,\ldots ,\Sigma_k$ be the sequence of all hexagons such that $(\Sigma_i)_t\cap a_{\pm n}\neq\emptyset$, for $i=1,2,\ldots ,k$. Note that $(\Sigma_1)_t\ni P_{\pm n}$ and $(\Sigma_k)_t\ni P_{\pm (n+1)}$. For a hexagon $\Sigma$, define $C(\Sigma )$ to be the union of all hexagons which are connected by a sequence of $0$-neighbors to $\Sigma$. Note that $C(\Sigma )$ looks like a trivalent tree and that it has infinitely many boundary components which are made out of partial boundaries of the hexagons in $C(\Sigma )$. It is important to note that either $C(\Sigma_1)=C(\Sigma_{k-1})$, or $C(\Sigma_1)$ and $C(\Sigma_{k-1})$ share a boundary component. If not,  then the subarc of $a_{\pm n}$ which connects $(\Sigma_1)_t$ to $(\Sigma_{k-1})_t$ connects two boundary components of some $C(\Sigma')$, where $C(\Sigma')$ separates $C(\Sigma_1)$ and $C(\Sigma_{k-1})$. 
Note that the arc which connects a short side of a hexagon to a non-adjacent side of the same hexagon has length at least $R/4$, where the long sides of the hexagon have lengths $R/2$. It follows that the subarc of  $a_{\pm n}$ which connects two boundary components of  $C(\Sigma')$ has length at least $R/4-3$.
Thus the above subarc of $a_{\pm n}$ has length greater than $1$ when $R$ is large enough which is impossible.

If $C(\Sigma_1)=C(\Sigma_{k-1})$ then we form a new sequence of hexagons $\Sigma_1,\Sigma_2',\ldots \Sigma_{k-2}',\Sigma_{k-1}$ such that the adjacent pairs of hexagons are $0$-neighbors  and $a_{\pm n}$ intersects characteristic triangles of the hexagons in the sequence. If $C(\Sigma_1)\neq C(\Sigma_{k-1})$ (and they share a boundary component) then we can choose a new sequence of hexagons $\Sigma_1,\Sigma_2',\ldots ,\Sigma_{k-2}',\Sigma_{k-1}$ such that each pair of adjacent hexagons are $0$-neighbors except one adjacent pair that are $1$-neighbors, and that $a_{\pm n}$ intersects characteristic triangles of the sequence. Note that the subarc of $a_{\pm n}$ that connects $(\Sigma_1)_t$ and $(\Sigma_{k-1})_t$ is of length less than $1$. 

The hexagons $\Sigma_{k-1}$ and $\Sigma_k$ are either $0$- or $1$-neighbors, or neither $0$- nor $1$-neighbors. If $\Sigma_{k-1}$ and 
$\Sigma_k$ are either $0$- or $1$-neighbors, then $\Sigma_1,\Sigma_2',\ldots ,\Sigma_{k-2}',\Sigma_{k-1},\Sigma_k$ is a sequence of hexagons whose adjacent hexagons are $0$-neighbors with the exception of at most $2$ pairs which are $1$-neighbors. Note that the arc $a_{\pm n}$ could have large length in general. If $\Sigma_k$ is a $0$-neighbor of $\Sigma_{k-1}$ then there is an arc $b_{\pm n}$ from the second point of the intersection of $a_{\pm n}$ with the boundary of $(\Sigma_{k-1})_t$ to the boundary of $(\Sigma_k)_t$ that has length less than $2$. To see this, let $\tilde{C}\in\pi^{-1}(\mathcal{P})=\tilde{\mathcal{P}}$ be the geodesic which contains one boundary side of both $\Sigma_{k-1}$ and $\Sigma_k$. Then the boundary side of $(\Sigma_{k-1})_t$ closets to $\tilde{C}$ is in the $C_1e^{-R/4}$-neighborhood of $\tilde{C}$ for some $C_1>0$, and the same statement is true for $(\Sigma_k)_t$. Since $\Sigma_k$ is shifted by $1\pm\frac{\epsilon}{R}$ with respect to $\Sigma_{k-1}$, it follows that such $b_{\pm n}$ exists. Thus the set of geodesics of $\tilde{\lambda}=\pi^{-1}(\lambda )$ that intersect $a_{\pm n}$ also intersect a geodesic arc $c_{\pm n}$ with the initial point $P_{\pm n}$ and of length at most $3$. Assume now that $\Sigma_{k}$ and $\Sigma_{k-1}$ are $1$-neighbors and that $\tilde{C}\in\tilde{P}$ separates them. Let $\Sigma_k'$ be the $0$-neighbor of $\Sigma_{k-1}$ which is separated by $\tilde{C}$ from $\Sigma_{k-1}$. It follows that the geodesics of $\tilde{\lambda}$ which intersect $a_{\pm n}$ except possibly the last geodesic (namely, the geodesic which contains one side of $(\Sigma_k)_t$ closets to $\tilde{C}$) intersect a geodesic arc of length at most $3$ with one endpoint $P_{\pm n}$. This follows simply by applying the above reasoning to the sequence $\Sigma_1,\Sigma_2',\ldots ,\Sigma_{k-2}',\Sigma_{k-1},\Sigma_k'$. If $\Sigma_k$ and $\Sigma_k'$ are not separated by some $\tilde{C}\in\mathcal{P}$, then $\Sigma_1,\Sigma_2',\ldots ,\Sigma_{k-2}',\Sigma_{k-1}$ suffices to get the same conclusion.

We give a proof of the nesting for the second case discussed above and the first case above is a subcase of the second. Namely, we are assuming that the set of geodesics $\tilde{\lambda}(a_{\pm n})$ of $\tilde{\lambda}$ which intersect $a_{\pm n}$ is also intersected by a geodesic arc $c_{\pm n}$ of length at most $3$ with the initial point $P_{\pm n}$ with a possible exception of one geodesic in $\tilde{\lambda}(a_{\pm n})$. We consider the bending map $(\tilde{f}_{\beta})|_{\Delta_{\Sigma_k}}=\varphi_{\Delta_{\Sigma_1},\Delta_{\Sigma_k}}$. Let $g_k$ be the geodesic of $\tilde{\lambda}$ which contains the boundary of $(\Sigma_k)_t$ and that separates $(\Sigma_k)_t$ and $(\Sigma_{k}')_t$. If $(\Sigma_k)_t\cap a_{\pm n}$ comes before $(\Sigma_{k}')_t\cap a_{\pm n}$ along $a_{\pm n}$ then  $c_{\pm n}$ does intersect $\Delta_{\Sigma_k}$ and this subcase of the second case reduces to the first case. Therefore, we assume that $(\Sigma_k)_t\cap a_{\pm n}$ comes after $(\Sigma_{k}')_t\cap a_{\pm n}$ along $a_{\pm n}$. The geodesic $g_k$ might not intersect $c_{\pm n}$. We have
$$
\varphi_{\Delta_{\Sigma_1},\Delta_{\Sigma_k}}=\varphi_{\Delta_{\Sigma_1},\Delta_{\Sigma_k'}}\circ R_{g_k}^{\beta (g_k)}
$$
where $\Sigma_k'$ is the $0$-neighbor of $\Sigma_{k-1}$ that is separated from $\Sigma_k$ by the geodesic $g_k$.

We normalize such that $P_{\pm n}=j\in\mathbb{H}^3$ and $P_{\pm (n+1)}=e^{-m}j$, where $m\geq 1$. Then $v=\{ e^{-m}j,-j\}$ is a tangent vector to $a_{\pm n}$ at the point $P_{\pm (n+1)}$ pointing towards the ideal endpoint of $g_{\pm 1}$. Lemma A.3 and the assumptions give
\begin{equation}
\label{eq:estimate_last_geodesic}
D_{T\mathbb{H}^3}(R_{g_k}^{\beta (\Delta_1(g_k),\Delta_2(g_k))}(\{ e^{-m}j,-j\}) ,\{e^{-m}j,-j\} )\leq C|\beta (\Delta_1(g_k),\Delta_2(g_k))|\leq \frac{C'\epsilon}{R}
\end{equation}
for some $C'>0$ when $\epsilon >0$ is small enough and $R\geq 1$, where $\{ e^{-m}j,-j\}\in T\mathbb{H}^3$ is a tangent vector to $\mathbb{H}^3$ based at $e^{-m}j$ and the function $D_{T\mathbb{H}^3}(\cdot ,\cdot )$ is defined in Appendix formula (\ref{eq:distance_on_H3}).

We consider $$\varphi_{\Delta_{\Sigma_1},\Delta_{\Sigma_{k}'}}=\psi_{\Delta_{\Sigma_1},\Delta_{\Sigma_k'}}\circ R_{g_k'}^{\beta (\Delta_{\Sigma_1},\Delta_{\Sigma_k'})}$$ where $g_k'\in\tilde{\lambda}$ is the side of $\Delta_{\Sigma_k'}$ facing $\Delta_{\Sigma_1}$. Let $\Delta_{\Sigma_i'}$ and $\Delta_{\Sigma_{i+1}'}$ be canonical triangles of two adjacent hexagons from the sequence $\Sigma_1':=\Sigma_1,\Sigma_2',\ldots ,\Sigma_{k-1}':=\Sigma_{k-1},\Sigma_k'$. Then  $\Delta_{\Sigma_i'}$ and $\Delta_{\Sigma_{i+1}'}$ are separated by $\tilde{C}_i\in\tilde{\mathcal{P}}$ and they have a common endpoint $\tilde{x}_i$ with $\tilde{C}_i$.
Let $\mathcal{F}_i$ be the family of complementary triangles between $\Delta_{\Sigma_i'}$ and $\Delta_{\Sigma_{i+1}'}$. Let $\Delta^i_1,\Delta^i_2$ be two triangles in $\mathcal{F}_i$ which are closets to $\Delta_{\Sigma_i'}$ and let $\Delta^{i+1}_1,\Delta^{i+1}_2$ be two triangles in $\mathcal{F}_i$ which are closets to $\Delta_{\Sigma_{i+1}'}$. Let $\gamma_i\in PSL_2(\mathbb{R})$ be the element of the covering group of $S$ that corresponds to $\tilde{C}_i$. Any triangle in $\mathcal{F}_i$
between $\Delta_{\Sigma_i'}$ and $\tilde{C}_i$ is the image of either $\Delta^i_1$ or $\Delta^i_2$ under a power of $\gamma_i$, and any triangle of $\mathcal{F}_i$ between $\tilde{C}_i$ and $\Delta_{\Sigma_{i+1}'}$ is the image of either $\Delta^{i+1}_1$ or $\Delta^{i+1}_2$ under a power of $\gamma_i$. 

\begin{figure}
\centering
\includegraphics[scale=0.5]{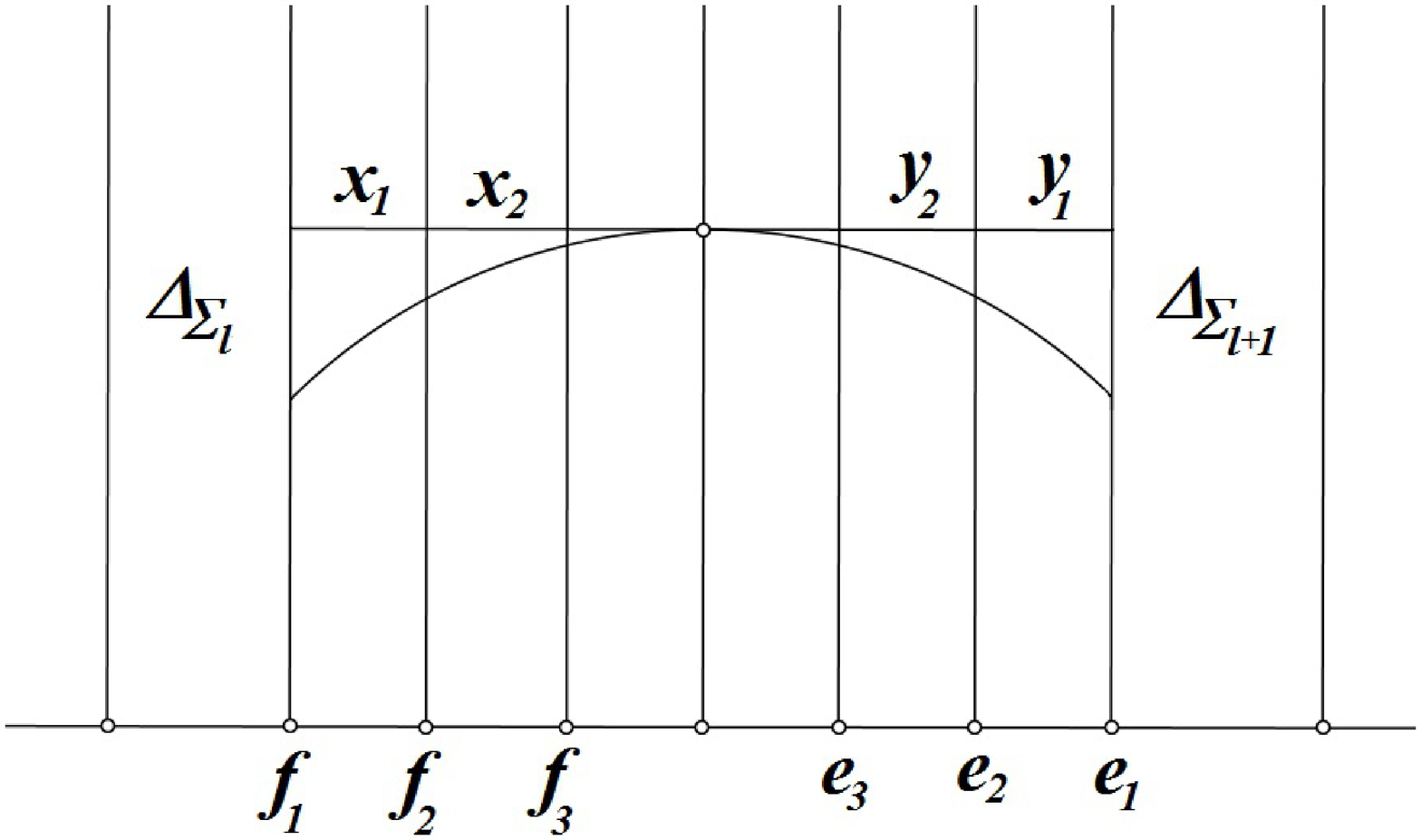}
\caption{}
\end{figure}

Let $h_i$ be the horocyclic arc connecting $\Delta_{\Sigma_i'}$ and $\Delta_{\Sigma_{i+1}'}$ with the center $\tilde{y}_i=c_{\pm n}\cap\tilde{C}_i$. Note that the length $|h_i|$ of the arc $h_i$ is less than a constant multiple of the length of the subarc of $c_{\pm n}$ connecting $\Delta_{\Sigma_i'}$ and $\Delta_{\Sigma_{i+1}'}$. Moreover, the length $h_i$ is less than the sum of the lengths of $h_i\cap\Delta_k^i$, $h_i\cap\Delta^{i+1}_k$, for $k=1,2$, and of the sum of the lengths of the intersections of $h_i$ with
the translates of $\Delta_k^i$ and $\Delta_k^{i+1}$, for $k=1,2$, under the powers of $\gamma_i$ such that the common endpoint of $\Delta^i_1$ and $\tilde{C}_i$ is repelling. Then the length of all translates is less than $C|h_i|e^{-l(\gamma_i)/2}$, where $l(\gamma_i)$ is the real part of the translation length of $\gamma_i$. To see this,
we normalize the
situation such that $\tilde{y}_i=(0,1)\in\mathbb{H}^2$ and $\tilde{C}_i$ is the
geodesic with endpoints $0$ and $\infty$. Then $h_i$ is the
horizontal Euclidean arc which contains
$(0,1)\in\mathbb{H}^2$ and the Euclidean length of
$h_i$ equals the hyperbolic length of $h_i$. Let $\{
f_m\}_{m\in\mathbb{N}}$ be the geodesics of
$\tilde{\lambda}$ with one endpoint $\infty$ that separate
$\Delta_{\Sigma_i'}$ and $\tilde{C}_i$ in the increasing order from
$\Delta_{\Sigma_i'}$. Let $\{ e_m\}_{m\in\mathbb{N}}$ be the
geodesics of $\tilde{\lambda}$ with one endpoint $\infty$
that separate $\Delta_{\Sigma_{i+1}'}$ and $\tilde{C}_i$ in the
decreasing order from $\Delta_{\Sigma_{i+1}}$. Let $x_1$ be the
length of the arc of $h_i$ between $f_1$ and $f_2$, and let $x_2$ be
the length of the arc of $h_i$ between $f_2$ and $f_3$. Note that
$f_m$ is mapped to $f_{m+2}$ by the hyperbolic translation $\gamma_i$ with the axis
$\tilde{C}_i$ and the attracting fixed point
$0\in\partial_{\infty}\mathbb{H}^2$. Then the distance between $f_{2m+1}$ and $f_{2m+2}$ is
$x_1e^{-\frac{ml(\gamma_i)}{2}}$, and similarly the distance between
$f_{2m+2}$ and $f_{2m+3}$ is $x_2e^{-\frac{ml(\gamma_i)}{2}}$ for
$m\in\mathbf{N}$ (cf. Figure 5). Therefore, the sum of the lengths of the gaps of $h_i$ except the first two gaps and the last two gaps is bounded by $C|h_i|e^{-l(\gamma_i)/2}$. 

Note that
\begin{equation}
\label{eq:bending_angle_finite}
\begin{array}l
\big{|}\beta (\Delta_{\Sigma_i'},\Delta^i_1)\big{|}\leq \frac{C_0\epsilon}{R}
\\
\big{|}\beta (\Delta_{\Sigma_i'},\Delta^i_2)\big{|}\leq \frac{2C_0\epsilon}{R}
\\
\big{|}\beta (\Delta_{\Sigma_i'},\Delta_{\Sigma_{i+1}'})\big{|}\leq \frac{C_0\epsilon}{R}
\\
\big{|}\beta (\Delta_{\Sigma_i'},\Delta_1^{i+1})\big{|}\leq \frac{2C_0\epsilon}{R}
\\
\big{|}\beta (\Delta_{\Sigma_i'},\Delta_2^{i+1})\big{|}\leq \frac{3C_0\epsilon}{R}.
\end{array}
\end{equation}

By the uniform boundedness of the composition of rotations \cite{Bon2}, there exists $C>0$ such that
$$
\|\psi_{\Delta_{\Sigma_1},\Delta_{\Sigma_k}}-id\|\leq C\sum_d\| R^{\beta (k_d)}_{g_d^{\Delta_{\Sigma_1}}} R^{-\beta (k_d)}_{g_d^{\Delta_{\Sigma_k}}}-id\|
$$
where the sum is over all gaps $d$ of $c_{\pm n}$, $k_d$ is the subarc of $c_{\pm n}$ from $P_{\pm n}$ to a point in $d$, and $g_d^{\Delta_{\Sigma_1}}$ ($g_d^{\Delta_{\Sigma_k}}$) is the leaf of $\tilde{\lambda}$ which contains the endpoint of $d$ closer to $\Delta_{\Sigma_1}$ ($\Delta_{\Sigma_k}$).  We divide the above sum over the gaps of $c_{\pm n}$ into two sums $\sum'$ and $\sum''$. The first sum $\sum'$ is over all gaps $c_{\pm n}\cap\Delta_{\Sigma_i'}$, for $i=1,2,\ldots ,k$, and $c_{\pm n}\cap \Delta^i_l$, for $l=1,2$, and the second sum $\sum''$ is over the remaining gaps. 

The first sum is finite. By Lemma \ref{lem:length1-intersects-at-most-R}, $k\leq 2R+2\leq 4R$ for $R\geq 1$ and by the finite additivity of $\beta$, we have that 
$$
\big{|}\beta (\Delta_{\Sigma_1},\Delta_{\Sigma_i'})\big{|}\leq C_1\epsilon
$$
for $i=1,2,\ldots ,k$, as well as
$$
\big{|}\beta (\Delta_{\Sigma_1},\Delta^i_l)\big{|}\leq C_1\epsilon
$$
for $i=1,2,\ldots ,k$ and $l=1,2$, and some constant $C_1>0$.

Lemma A.4 implies
$$
\sum'\leq \sum_{i=1}^kC_2\epsilon |h_i|\leq C_2 |c_{\pm n}|\epsilon\leq C_3\epsilon .
$$
It remains to estimate $\sum''$. We proved above that the total length of the gaps of $h_i$ with respect to the family $\mathcal{F}_i$ of complementary triangles except for the first two and the last two gaps is less than $Ce^{-R/2}|h_i|$. Since $\sum_{i=1}^{n}|h_i|\leq C|c_{\pm n}|\leq C_4$ and the cocycle $\beta$ takes values in $[-\pi ,\pi)$ (thus $\beta$ is bounded),
Lemma A.4 gives
$$
\sum''\leq C\sum_{i=1}^ke^{-R/2}\leq C_5Re^{-R/2}
$$
for some $C_5>0$. Then $\sum' +\sum''$ can be made arbitrary small when $\epsilon >0$ is small enough and $R>0$ is large enough. The above, Lemma A.3 and 
$$
\big{|}\beta (\Delta_{\Sigma_1},\Delta_{\Sigma_k'})\big{|}\leq C\epsilon
$$
imply that 
$$
\|\varphi_{\Delta_{\Sigma_1},\Delta_{\Sigma_k'}}-id\|
$$
is as small as needed for $\epsilon >0$ small enough and $R>0$ large enough. Then (\ref{eq:bending_angle_finite}), Lemma A.2 and the above prove that the assumptions of Lemma A.1 are satisfied for $\epsilon >0$ small enough and $R>0$ large enough. Thus the nesting for $\tilde{f}_{\beta}$ on $a_{\pm n}$ follows by Lemma A.1. We choose $\hat{\epsilon}>0$ and $R(\hat{\epsilon})>0$ accordingly.

We assume now that $\Sigma_{k-1}'$ and $\Sigma_k$ are neither $0$- nor $1$-neighbors. Then there is a unique $\tilde{C}_{k-1}\in\tilde{\mathcal{P}}$ which separates $\Sigma_{k-1}'$ and $\Sigma_k$, and that contains boundary sides of both of them. Let $\Sigma_k'$ be the $0$-neighbor of $\Sigma_{k-1}'$ which is on the same side of $\tilde{C}_{k-1}$ as $\Sigma_k$. Let $s\geq 1$ be the number of hexagons  in between $\Sigma_k'$ and $\Sigma_k$. There are two possibilities: either $a_{\pm n}$ intersect $\Delta_{\Sigma_k'}$ in which case we say that $\Delta_{\Sigma_k}$ is ``above'' $\Delta_{\Sigma_k'}$, or $a_{\pm n}$ does not intersect $\Delta_{\Sigma_k'}$ in which case we say that $\Delta_{\Sigma_k}$ is ``below'' $\Delta_{\Sigma_{k-1}'}$. 

Assume we are in the former case and let $\{ g_1,g_2,\dots ,g_{s+1}\}$ be the geodesic of $\tilde{\lambda}$ between $\Delta_{\Sigma_k'}$ and $\Delta_{\Sigma_k}$. We use the following fact. Let $h$ and $h'$ be two geodesics that intersect $L=\{ (0,0,t):t>0\}$ at points $e^{-m}j$ and $e^{-m'}j$ subtending angles $\epsilon >0$ and $\epsilon'>0$, where $m<m'$. Let $\epsilon''=\max\{\epsilon ,\epsilon'\}$ and let $h''$ be the geodesic that intersects $L$ at the point $e^{-m}j$ subtending an angle $\epsilon''$. Then, for $m''\geq m'$ and $\epsilon''>0$, we have
\begin{equation*}
\begin{split}
D_{T\mathbb{H}^3}(R_h^{\theta}\circ R_{h'}^{\theta'}(\{ e^{-m''}j,-j\}), \{ e^{-m''}j,-j\})\leq\ \ \ \ \ \ \ \ \  \\ 
\leq\max_{0\leq\theta''\leq 2\pi}D_{T\mathbb{H}^3}(R_{h''}^{\theta''}(\{ e^{-m''}j,-j\}),\{ e^{-m''}j,-j\}).
\end{split}
\end{equation*}

Let $$R_s=R_{g_1}^{\beta (\Delta_1(g_1),\Delta_2(g_1))}\circ\cdots\circ R_{g_{s+1}}^{\beta (\Delta_1(g_{s+1}),\Delta_2(g_{s+1}))}.$$ 
Then, the above implies that
\begin{equation}
\label{eq:dist_max}
\begin{split}
D_{T\mathbb{H}^3}(R_s(\{ e^{-m''}j,-j\}), \{ e^{-m''}j,-j\})\leq\ \ \ \ \ \ \ \ \ \ \ \ \ \ \ \ \\
\leq\max_{0\leq\theta\leq 2\pi}D_{T\mathbb{H}^3}(R_{g_1'}^{\theta}(\{ e^{-m''}j,-j\}),\{ e^{-m''}j,-j\}).
\end{split}
\end{equation}
where $g_1'$ is the geodesic passing through $g_1\cap L$ that subtends an angle $\max\{ |\angle (g_1,L)|,\ldots ,|\angle (g_s,L)|\}$ with $L$.
Lemma A.5 and (\ref{eq:dist_max}) imply that
$D_{T\mathbb{H}^3}(R_s(\{ e^{-m''}j,-j\}), \{ e^{-m''}j,-j\})$ is as small as we want when the angle $|\angle (g_1',L)|$ is small enough for any $0\leq\theta\leq 2\pi$.
Note that 
$$
\varphi_{\Delta_{\Sigma_1},\Delta_{\Sigma_k}}=\varphi_{\Delta_{\Sigma_1},\Delta_{\Sigma_k'}}\circ \varphi_{\Delta_{\Sigma_k'},\Delta_{\Sigma_k}}.
$$
Since $\varphi_{\Delta_{\Sigma_k'},\Delta_{\Sigma_k}}=R_s$, the above gives
$$
D_{T\mathbb{H}^3}(\varphi_{\Delta_{\Sigma_k'},\Delta_{\Sigma_k}}(\{ P_{\pm (n+1)},-j\} ),\{ P_{\pm (n+1)},-j\})
$$
is as small as needed for $R$ large enough. Indeed, the subarc of $a_{\pm n}$ from the second point of the intersection with the boundary of $(\Sigma_{k-1})_t$ to the first point of intersection with the boundary of $(\Sigma_k')_t$ is inside one complement of ${\mathcal TH}_t$ as well as long sub arcs of the set of geodesics $\{ g_1,\ldots ,g_{s+1}\}$. Thus $a_{\pm n}$ and $\{ g_1,\ldots ,g_{s+1}\}$ remain in a neighborhood of one $\tilde{C}\in\tilde{\mathcal{P}}$ for a long distance when $R$ is large. It follows that the angles of intersections between $a_{\pm n}$ and the geodesics in $\{ g_1,g_2,\ldots ,g_{s+1}\}$ are small for $R$ large enough and the above applies. The reasoning in the first case applies to $\varphi_{\Delta_{\Sigma_1},\Delta_{\Sigma_k'}}$ and we have the nesting of the images of the cones at the endpoints of $a_{\pm n}$ under the bending map $\varphi_{\Delta_{\Sigma_1},\Delta_{\Sigma_k}}$.
If $\Delta_{\Sigma_{k-1}'}$ is ``above'' $\Delta_{\Sigma_k}$ then symmetry reduces to the previous case.

It remains to consider the case when $a_{\pm n}$ has infinite length (in which case the endpoint of $a_{\pm n}$ is also the endpoint of $\tilde{C}\in\tilde{\mathcal{P}}$ and $a_{\pm n}\subset\mathbb{H}^2-\mathcal{TH}_t$).  
An elementary (euclidean) considerations prove that when $R\geq 1$ the number of geodesics of $\tilde{\lambda}$ that intersect the geodesics subrays of $a_{\pm n}$ which connect two sides of a single hexagon is at most $6$. Indeed, assume that $a_{\pm n}$ is the geodesic arc in $\mathbb{H}^2$ with the initial point $i-e^{-R/4}$ and the endpoint $0\in\partial_{\infty}\mathbb{H}^2$. Then $a_{\pm n}$ is a circular arc with the center $a=\frac{e^{R/4}+e^{-R/4}-\sqrt{(e^{R/4}+e^{-R/4})^2-4e^{-R}}}{2}$ and the radius $a$. The $x$-coordinate of the intersection of $a_{\pm n}$ with the horizontal line $y=e^{-R/2}$ is estimated to be more than $e^{-3R}{2}$. Since the translation length of the element $\gamma$ fixing the $y$-axis is $e^{-R/2}$ and since $\gamma$ identifies every second geodesic of $\tilde{\lambda}$ that have endpoint $\infty$, the claim follows.
Then applying Lemma A.3 finitely many times to the sequence of subarcs of $a_{\pm n}$ of lengths $R/2$, we obtain a nesting property along this sequence. Thus $\tilde{f}_{\beta}$ is injective for $\epsilon >0$ small enough and $R>0$ large enough. 

We choose $\hat{\epsilon}$ and $R(\hat{\epsilon})$ as the minimum of the choice in all the cases considered.
\end{proof}

\section{Holomorphic motions}

We finish the proof of Theorem \ref{thm:injective-bending-cocycles} using holomorphic motions. This proof is standard once the injectivity is established (cf. \cite{KeeSe}, \cite{EpMarMar} and \cite{Sa1}). 
Holomorphic motions were introduced and studied in \cite{MSS} and the key extension property is proved in \cite{Slo}. 

The endpoints of the representations of elements of $\pi_1(S)$ vary holomorphically in the complex Fenchel-Nielsen coordinates. We established that the holomorphic variation is injective on the set of endpoints when the parameters are close to being real in the sense of (\ref{eq:twist-bend}) and (\ref{eq:half-length-modif}). Thus the holomorphic variation of the endpoint of $\pi_1(S)$ is a holomorphic motion which extends by the lambda lemma (cf. \cite{MSS}) to a holomorphic motion of the unit circle. Then there exists an extension to a holomorphic motion of the complex plane (cf. \cite{Slo}). It follows that $\tilde{f}_{\beta}$ extends to a quasiconformal mapping of the complex plane 
and that the quasiconformal constant is less than $1+K_0\epsilon$ for $\epsilon >0$ small enough and fixed $K_0>0$ (cf. \cite{MSS}). The extension of $\tilde{f}_{\beta}$ can be chosen to be equivariant with respect to to the action of $\pi_1(S)$ (cf. \cite{EKK}) which finishes the proof of Theorem \ref{thm:injective-bending-cocycles}. 

\section*{Appendix}

We use the quaternions to represent the upper half-space model
$\mathbb{H}^3=\{ z+tj:z\in\mathbb{C},t>0\}$ of the hyperbolic
three-space (see Beardon \cite{Bear}), where $j=(0,0,1)\in\mathbb{H}^3$. The space of isometries of
$\mathbb{H}^3$ is identified with $PSL_2(\mathbb{C})$.
The Poincar\'e extension of $A(z)=\frac{az+b}{cz+d}\in PSL_2(\mathbb{C})$ to $\mathbb{H}^3$  is given in \cite{Bear}
by
$$A(z+tj)=\frac{(az+b)\overline{(cz+d)}+a\bar{c}t^2+tj}{|cz+d|^2+|c|^2t^2}.
$$

An isometry of $\mathbb{H}^3$ which is close to the identity moves
points on a bounded distance from $j\in\mathbb{H}^3$ by a small
amount and the tangent vectors are rotated by a small angle with
respect to the Euclidean parallel transport in $\mathbb{R}^3$. We
give a quantitative statement for the above including the situation
when the points are on the unbounded distances from
$j\in\mathbb{H}^3$ which is needed in our considerations.

Given $P=z+tj\in\mathbb{H}^3$, we define
$$
ht(P)=t
$$
and
$$
Z(P)=z.
$$

Consider the tangent space $T\mathbb{H}^3$ to the upper half-space $\mathbb{H}^3$. Let $\{ P,u\} ,\{ Q,v\}\in T\mathbb{H}^3$ be two tangent vectors based at $P,Q\in\mathbb{H}^3$, respectively. 
We define
\begin{equation}
\label{eq:distance_on_H3}
D_{T\mathbb{H}^3}(\{ P,u\} ,\{ Q,v\} )=\max\{ \big{|}\frac{ht(P)}{ht(Q)}-1\big{|}, \big{|}Z(P)-Z(Q)\big{|},\big{|}\angle (u,v)\big{|}\},
\end{equation}
where $\angle (u,v)$ is the angle between the vectors $u$ and $v$ after the euclidean transport in $\mathbb{H}^3$. Note that  $D_{T\mathbb{H}^3}(\{ P,u\} ,\{ Q,v\} )$ is not a metric on $T\mathbb{H}^3$. 

\vskip .2 cm

\paragraph{\bf Lemma A.1}
{\it 
Given $m_0>0$, there exists $\delta =\delta (m_0)>0$ such that for any $m\geq m_0$ we have
\begin{equation}
\label{eq:nested-cone}
\overline{\partial_{\infty}\mathcal{C}(P,v,\frac{\pi}{2})}\subset
\partial_{\infty}\mathcal{C}(j,-j,\frac{\pi}{2})
\end{equation}
where $\{ P,v\}\in T\mathbb{H}^3$ satisfies $D_{T\mathbb{H}^3}(\{ P,v\} ,\{ e^{-m}j,-j\} )< \delta$.
}

\vskip .2 cm

\begin{proof}
Let $x$ be the center of the Euclidean hemisphere (orthogonal to $\mathbb{C}$) that passes through $P$ and touches the unit Euclidean hemisphere orthogonal to $\mathbb{C}$ with the center $0\in\mathbb{C}\subset\partial_{\infty}\mathbb{H}^3$ (cf. Figure 6).
Let $y=Z(P)\in\mathbb{C}$ and let $\varphi$ be the angle between euclidean segments $Px$ and $Py$ at the point $P$.

An elementary (Euclidean) considerations give 
\begin{equation}
\label{eq:x}
x\geq C_1(m_0)>0
\end{equation}
for some constant $C_1(m_0)>0$ which depends on $m_0$.

\begin{figure}
\centering
\includegraphics[scale=0.5]{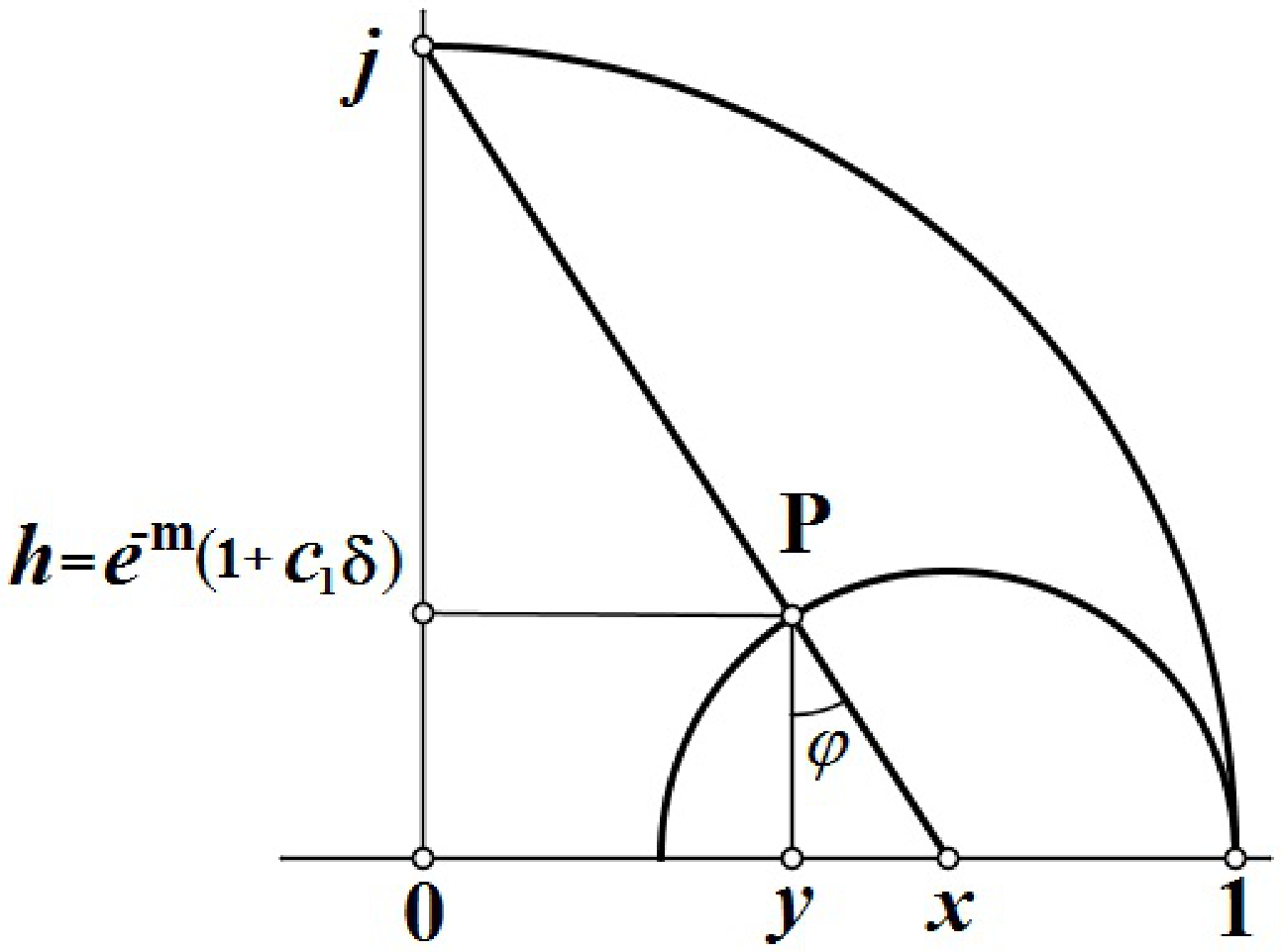}
\caption{}
\end{figure}

This implies that 
$$
\varphi \geq C_2(m_0)>0
$$
for some constant $C_2(m_0)>0$.
Thus
\begin{equation}
\label{eq:inter_nesting}
\overline{\partial_{\infty}\mathcal{C}(P,-j,\frac{\pi}{2}+C_2(m_0))}\subset \partial_{\infty}\mathcal{C}(j,-j,\frac{\pi}{2}).
\end{equation}
Since the angle (at the point $P$) between the hyperbolic geodesic connecting $j$ to $P$ and the euclidean segment $Py$ is less than $\varphi$, the above inclusion implies (\ref{eq:nested-cone}) for $\delta (m_0)<C_2(m_0)$.
\end{proof}

\vskip .2 cm

\paragraph{\bf Lemma A.2}
{\it Let $g\in PSL_2(\mathbb{C})$ with
$$
\| g-id\| <\frac{1}{2}
$$
and let $\{ z+tj,u\}$ be a tangent vector to $\mathbb{H}^3$ such that
$$
D_{T\mathbb{H}^3}(\{ z+tj,u\} ,\{ e^{-m}j,-j\} )\leq\delta
$$
for $0\leq\delta <\delta_0$ with $\delta_0 >0$ fixed.
Then there exist $C_1,C_2>0$ depending on $\delta_0$ such that
for every $m>0$ we have
$$
D_{T\mathbb{H}^3}(g(\{ z+tj,u\}),\{ e^{-m}j,-j\} )\leq C_1\delta +C_2\|g-id\|.
$$}

\vskip .2 cm

\begin{proof} 
Denote by $g$ the Poincar\'e extension of $g(z)=\frac{az+b}{cz+d}$ with $ad-bc=1$ given above. Then
$$
\Big{|}ht(g(z+tj))-e^{-m}\Big{|}=\Big{|}
\frac{t}{|cz+d|^2+|c|^2t^2}-e^{-m}\Big{|}\leq
C_1e^{-m}\delta +C_2e^{-m}\| g-id\|
$$
for constants $C_1,C_2>0$ independent of $m>0$ and for all $g\in
PSL_2(\mathbf{C})$ with $\| g-id\|\leq\frac{1}{2}$. Moreover,
$$
\Big{|}Z(g(z+tj))\Big{|}=\Big{|}\frac{a\bar{c}|z|^2+a\bar{d}z+b\bar{c}\bar{z}+b\bar{d}+a\bar{c}t^2}{|cz+d|^2+|c|^2t^2}\Big{|}
\leq C_1\delta +C_2\| g-id\|.
$$

Let $u=<u_1,u_2,u_3>$. Without loss of generality, we assume that $|u_1|,|u_2|\leq\delta$ and $u_3=-j$. Let $v=Dg(z+tj)u=<v_1,v_2,v_3>$. Let $g=g_1+g_2i+g_3j$ be the
coordinate functions of $g$. Direct computations give
$$
\big{|}\frac{\partial g_i}{\partial x}(z+tj)\big{|} ,\big{|}\frac{\partial g_i}{\partial y}(z+tj)\big{|}\leq C
$$
for some $C>0$ and $i=1,2,3$, where $z+tj=x+yi+tj$. Moreover,
$$
\big{|}\frac{\partial g_i}{\partial t}(z+tj)\big{|}\leq C_1\| g-id\|
$$
for some $C_1>0$ and $i=1,2$, and
$$
\big{|}\frac{\partial g_3}{\partial t} (z+tj)\big{|}\geq 1- C_2\| g-id\|
$$
for some $C_2>0$. 

The above inequalities give the following
$$
|v_1|,|v_2|\leq C'(\delta +\| g-id\|)
$$
and $$v_3\leq -1+C''(\delta +\| g-id\|).$$
This gives that 
$$
\big{|}\angle (-j,v)\big{|}\leq C'''(\delta +\| g-id\|).
$$
The lemma is proved.
\end{proof}

Let $L=\{ (0,0,t):t>0\}\subset\mathbb{H}^3$ be the geodesic through $j=(0,0,1)\in\mathbb{H}^3$ with the ideal endpoint $0\in\mathbb{C}\subset\partial_{\infty}\mathbb{H}^3$.

\vskip .2 cm

\paragraph{\bf Lemma A.3}
{\it
Let $h$ be a geodesic in $\mathbb{H}^2\subset\mathbb{H}^3$ that intersects $L$ between points $j$ and $e^{-r}j$ for some $r\geq 1$. Given $\epsilon_0,\delta_0>0$, there exist $C(r,\epsilon_0,\delta_0), C_0(r,\epsilon_0,\delta_0) >0$
such that 
$$
D_{T\mathbb{H}^3}(R_h^{\epsilon}(\{ z,u\} ),\{e^{-r}j,-j\})\leq C_0\delta +C\epsilon
$$
for any $0\leq\epsilon <\epsilon_0$, $0\leq\delta <\delta_0$, and $\{ z,u\}\in T\mathbb{H}^3$ with
$$D_{T\mathbb{H}^3}(\{ z,u\} ,\{e^{-r}j,-j\} )\leq\delta ,$$
where $R_h^{\epsilon}$ is the hyperbolic rotation around the axis $h$ by the angle $\epsilon$.
}

\begin{proof}
The quantity $D_{T\mathbb{H}^3}(R_h^{\epsilon}(\{ z,u\} ),\{e^{-r}j,-j\})$ is the largest when $h$ is orthogonal to $L$ at the point $j$. In this case $R_h^{\epsilon}$ fixes $1$ and $-1$, and
$$
R_h^{\epsilon}(z)=\frac{\cos\frac{\epsilon}{2}z-i\sin\frac{\epsilon}{2}}{-i\sin\frac{\epsilon}{2}z+\cos\frac{\epsilon}{2}}.
$$
Therefore, there exists $C>0$ such that
$$
\| R_h^{\epsilon}-id\|\leq C\epsilon.
$$
The lemma follows by Lemma A.2.
\end{proof}

\vskip .2 cm

The following lemma is standard (cf. \cite{Bon2}).

\vskip .2 cm

\paragraph{\bf Lemma A.4}
{\it
Let $D_{r_0}(j)\subset\mathbb{H}^3$ be the hyperbolic ball of radius $r_0>0$ centered at $j\in\mathbb{H}^3$. Then there exists $C=C(r_0)>0$ such that if $h_1,h_2$ are two hyperbolic geodesics with a common endpoint that intersect $D_{r_0}(j)$ and if $d_{r_0}(h_1,h_2)$ is the hyperbolic distance between $h_1\cap D_{r_0}(j)$ and $h_2\cap D_{r_0}(j)$ then
$$
\| R^{\epsilon}_{h_1}R^{-\epsilon}_{h_2}-id\| \leq Cd_{r_0}(h_1,h_2)\epsilon
$$
for any $\epsilon >0$. 
}

\vskip .2 cm

\paragraph{\bf Lemma A.5}
{\it
Let $g$ be a geodesic in $\mathbb{H}^2\subset\mathbb{H}^3$ that intersects $L=\{ (0,0,t):t>0\}$ between $j$ and $e^{-m}j$ at an angle $\epsilon >0$. Then, for any $\delta >0$,
$$
D_{T\mathbb{H}^3}(\{ R_g^{\theta}(\{ e^{-m}j,-j\}),\{ e^{-m}j,-j\} )<\delta
$$
when $\epsilon =\epsilon (\delta )>0$ is small enough.
}

\vskip .2 cm

\begin{proof}
The angle $\epsilon_1>0$ between $R_g^{\theta}(L)$ and $L$ satisfies
\begin{equation}
\label{eq:angle_intersection_small}
\big{|}\epsilon_1\big{|}\leq\frac{\pi}{2}\cdot\big{|}\theta\big{|}\cdot\big{|}\epsilon\big{|}.
\end{equation}
To see this, we normalize such that $g=L$ and $L$ has turned into a geodesic with endpoints $a<0$ and $b>0$. The the geodesic $L$ is parametrized by $\gamma (t)=\frac{a+b}{2}+\frac{b-a}{2}\cos t+(\frac{a+b}{2}+\frac{b-a}{2}\sin t)j$ and $R_g^{\theta}(L)$ is parametrized by $\gamma_1 (t)=(\frac{a+b}{2}+\frac{b-a}{2}\cos t)\cos\theta+(\frac{a+b}{2}+\frac{b-a}{2}\cos t)\sin\theta i +(\frac{a+b}{2}+\frac{b-a}{2}\sin t)j$. Thus 
$$
\cos\epsilon_1 =\frac{\gamma'(t)\cdot\gamma_1'(t)}{\|\gamma'(t)\|\cdot\|\gamma_1'(t)\|}=\sin^2\epsilon\cos\theta +\cos^2\epsilon
$$
which implies
$$
|\epsilon_1|\leq |2\sin^{-1}(\epsilon\sin\frac{\theta}{2})|\leq\frac{\pi}{2}|\theta |\cdot |\epsilon |.
$$

We go back to the assumption that $L=\{ (0,0,t):t>0\}$. By (\ref{eq:angle_intersection_small}), the absolute value of the angle $|\epsilon_1|$ between $L$ and $R_g^{\theta}(L)$ is less than $\frac{\pi^2}{2}|\epsilon |$. Let $P=g\cap L$. The angle (after the Euclidean parallel transport) between the tangent vector to $R_g^{\theta}([P,e^{-m'}j])$ at the endpoint $e^{-m'}j$ and the vector $-j$ decreases as $m'\to m$ for all $m'$ such that $ht(P)\geq ht(e^{-m'}j)\geq ht(e^{-m}j)$. Thus the angle is the largest if $P=e^{-m}j$ which implies that $|\epsilon_1|\leq\frac{\pi^2}{2}|\epsilon |$.

Thus, for $R$ large enough, we have
$$
\big{|}\angle (R_g^{\theta}(\{ e^{-m}j,-j\}),\{ e^{-m}j,-j\})\big{|}<\delta .
$$  
To estimate the height and $z$-coordinate of $R_g^{\theta}(e^{-m}j)$, we note that both
$|ht(R_h^{\theta}(e^{-m}j))-ht(e^{-m}j)|$ and
$|Z(R_h^{\theta}(e^{-m}j))|$ are the largest when $g\cap h=p=j$.
Then the angle between $g$ and $R_h^{\theta}(g)$ is $\epsilon_1$,
where $|\epsilon_1|\leq \frac{\pi^2}{2}|\epsilon |$. Since $g$ and
$R_h^{\theta}(g)$ belong to a hyperbolic plane embedded in
$\mathbb{H}^3$ which contains the $g$, we can restrict further
analysis to the upper half-plane $\mathbb{H}^2$ where we identify
$i\in\mathbb{H}^2$ with $j\in\mathbb{H}^3$. Let $A\in
PSL_2(\mathbb{R})$ be an isometry of $\mathbb{H}^2$ which fixes
$i\in\mathbb{H}^2$ and maps $g$ into $R_h^{\theta}(g)$. Then
$A(e^{-m}i)=R_h^{\theta}(e^{-m}j)$ for the embedding
$\mathbb{H}^2\subset\mathbb{H}^3$. Note that
$A(z)=\frac{cz+d}{-dz+c}$ with $c,d>0$ and $c^2+d^2=1$. It follows
that $|\frac{d}{c}|\leq C\epsilon_1$, for some $C>0$ and for $|\epsilon |$
small enough. Furthermore,
$$
A(e^{-m}i)=\frac{e^{-m}i+\frac{d}{c}}{-\frac{d}{c}e^{-m}i+1}
$$
which implies that
$$
\Big{|}Im(A(e^{-m}i))-e^{-m}\Big{|}\leq
C_{11}\Big{(}\frac{d}{c}\Big{)}^2e^{-m}
$$
and
$$
\Big{|}Re(A(e^{-m}i))\Big{|}\leq C_{12}\frac{d}{c}
$$
and the lemma follows.
\end{proof}

\end{document}